\documentclass[a4paper,11pt,reqno]{article}
\usepackage{amsmath}
\usepackage{amssymb,latexsym}
\usepackage[english]{babel}
\usepackage[latin1]{inputenc} 
\usepackage[T1]{fontenc} 
\usepackage{ae}
\usepackage{color,graphicx}
\usepackage{subfigure}

\newcommand{\goodgap}{\hspace{\subfigtopskip}\hspace{\subfigbottomskip}}

\newdimen\AAdi%
\newbox\AAbo%

\def\AArm{\fam0 \rm}%
\def\AAk#1#2{\setbox\AAbo=\hbox{#2}\AAdi=\wd\AAbo\kern#1\AAdi{}}%
\def\AAr#1#2#3{\setbox\AAbo=\hbox{#2}\AAdi=\ht\AAbo\raise#1\AAdi\hbox{#3}}%

\newcommand{\un}{{\AArm 1\AAk{-.9}{l}l}}
\renewcommand{\Re}{\mathrm{Re\:}}

\newtheorem{thm}{Theorem}[section]

\newtheorem{lem}[thm]{Lemma}
\newtheorem{prop}[thm]{Proposition}
\newtheorem{rem}[thm]{Remark}

\def\qedbox{$\rlap{$\sqcap$}\sqcup$}
\def\finpreuve{\nobreak\hfill\penalty250 \hbox{}
           \nobreak\hfill\qedbox\par\medskip}

\newcommand{\J}{J_t\,,\;t\,{\ge}\,0}

\newcommand{\B}{B_t\,,\;t\,{\ge}\,0}
\newcommand{\X}{X_t\,,\;t\,{\ge}\,0}

\newcommand{\Btild}{\widetilde B_t\,,\;t\,{\ge}\,0}

\def\sectio#1{\setcounter{equation}{0}\section{#1}}
\def\eqp#1{\ensuremath\,#1} 

\title{\rm \Large \bf Lévy Processes:\\ Hitting Time,
Overshoot and Undershoot\\
I -  Functional Equations} \date{Manuscript no.\ P660 submitted to SPA, October 2004}

\author{Bernard Roynette$^a$ - Pierre Vallois$^a$ - Agn\`es Volpi$^b$}

\begin{document}

\maketitle
\vspace{-0.5cm}

\noindent $^a$D\'epartement de math\'ematique, Institut {\'E}lie
Cartan,Universit{\'e} Henri Poincar{\'e}, BP 239, 54506 Vand\oe uvre-l\`es-Nancy cedex, France.

\noindent $^b$ESSTIN, 2 rue Jean Lamour, Parc Robert Bentz, 54500 Vand\oe uvre-l\`es-Nancy, France.

\noindent e-mail:\\
roynette@iecn.u-nancy.fr~~vallois@iecn.u-nancy.fr~~agnes.volpi@esstin.uhp-nancy.fr

\vspace{1ex}

\begin{abstract}

Let ($X_t, \; t\ge0$) be a Lévy process started at $0$, with Lévy measure $\nu$, and $T_x$ the first hitting time of level $x>0$~: 
\mbox{$\displaystyle T_x:=\inf{ \{ t\ge 0;\;  X_t>x \}}$}. Let $F(\theta,\mu,\rho,.)$ be the joint Laplace transform of $(T_x, K_x, L_x)$~: \goodbreak \noindent
$\displaystyle 
F(\theta,\mu,\rho,x):=\mathbb E \left(e^{-\theta T_x-\mu K_x-\rho L_x}\un_{\{T_x<+\infty\}}\right) \eqp{,}
$
where $\theta\geq 0$, $\mu\geq 0$, $\rho\geq 0$, $x\geq 0$, 
  \mbox{$\displaystyle K_x:=X_{T_x}-x$}  and  \mbox{$\displaystyle L_x:=x-X_{T_{x^-}}$}.\\
If $\nu(\mathbb R)<+\infty$ and $\displaystyle \int_1^{+\infty} e^{sy}\nu(dy)<+\infty$ for some $s>0$, then we prove that $F(\theta,\mu,\rho,.)$ is the unique solution of an integral equation and has a subexponential decay at infinity when $\theta>0$ or $\theta=0$ and $\mathbb E(X_1)<0$. If $\nu$ is not necessarily a finite measure but verifies $\displaystyle \int_{-\infty}^{-1} e^{-sy}\nu(dy)<+\infty$ for any $s>0$, then the $x$-Laplace transform of $F(\theta,\mu,\rho,.)$ satisfies some kind of integral equation. This allows us to prove that $F(\theta,\mu,\rho,.)$ is a solution to a second integral equation.

\end{abstract}

\vspace{1ex}



\vspace{1ex}

\section*{Introduction}
\addcontentsline{toc}{section}{Introduction}

\begin{enumerate}
\item  Let ($X_t, \; t\ge0$) be a Lévy process, right continuous with left limits,  started at $0$, with Lévy measure $\nu$. 
We suppose that ($\X$) may  be decomposed as follows:
\begin{equation}
\label{intro:def-de-X}
 X_t=\sigma B_t-c_0t+J_t \qquad t\geq 0\eqp{,}
\end{equation}
where $c_0\in \mathbb R$, $\sigma >0$,  ($\B$) is an one-dimensional Brownian motion started at $0$,  ($\J$) is a pure jump Lévy process, ind\'ependent of ($\B$) and $J_0=0$. We will suppose that $\sigma=1$.\\
Let us introduce the function $\varphi$ which will play a central role in our study~: $\varphi(q):=\psi(-q)$ where $\psi$ is the characteristic exponent of $(\X)$, i.e. \mbox{$\displaystyle \mathbb E(e^{q X_t})=e^{t\psi(q)}$}. By Lévy-Khintchine formula, we have~:
\begin{equation}
\label{def-phi}
\varphi(q)=\frac{q^2}{2}+cq+\int_{\mathbb R}(e^{-qy}-1+qy \un_{\{|y|<1\}})\nu(dy)\eqp{.}
\end{equation}
Remark if~:
\begin{equation}
\label{moment1-nu}
\int_\mathbb R |y|\un_{\{|y|>1\}}\nu(dy)<+\infty\eqp{,}
\end{equation}
then $X_1$ has a finite expectation and~:
\begin{equation}
\mathbb E(X_1)=-c+\int_\mathbb R y\un_{\{|y|\geq 1\}}\nu(dy)\eqp{.}
\end{equation}
\item In this paper we are interested in the first hitting time  of level $x>0$ :
\begin{equation}
 T_x:=\inf{ \{ t\ge 0;\;  X_t>x \}}\eqp{.}
\end{equation}

Setting $Z_t:=x-X_t$, then  
\mbox{$\displaystyle T_x:=\inf{ \{ t\ge 0;\;  Z_t<0 \}}$} is the ruin time to a company whose fortune is  modelled by $(Z_t\,;\;t \geq 0)$. 

We also consider the overshoot $K_x$, respectively the undershoot $L_x$ :
\begin{align}
&\label{def-overshoot}
 K_x:=X_{T_x}-x \eqp{,}\\
&\label{def-undershoot}
L_x:=x-X_{T_{x^-}}\eqp{.}
\end{align}

\noindent The aim of this paper is the study of the joint distribution of  $(T_x, K_x, L_x)$.

\noindent Our approach makes appeal to the joint Laplace transform of $(T_x, K_x, L_x)$, namely, for all $\theta\geq0$, $\mu\geq0$, $\rho\geq0$, $x\geq0$~: 
\begin{equation}
\label{def-de-F-theta-mu-rho}
F(\theta,\mu,\rho,x)=\mathbb E \left(e^{-\theta T_x-\mu K_x-\rho L_x}\un_{\{T_x<+\infty\}}\right) \eqp{.}
\end{equation}
If $\theta=\mu=\rho=0$, 
\begin{equation}
\label{def-de-proba-de-ruine}
F(0,0,0,x)=\mathbb P(T_x<+\infty) 
\end{equation}
is the well-known  ruin probability.\\

\item In section~\ref{sec:equa-fonctionnelle} we study $F(\theta,\mu,\rho,.)$ when $\nu(\mathbb R)<+\infty$ (i.e. ($\J$) is a compound Poisson process). Since ($\X$) has a first jump time $\tau_1$, and $(X_{t+\tau_1}-X_{\tau_1}\;;\;t\geq0)$ is distributed as ($\X$), we show  in Theorem~\ref{thm:equa-fonctionnelle} that $F(\theta,\mu,\rho,.)$ verifies the integral equation (\ref{equa-fonctionnelle}). To go further we suppose moreover~:
\begin{equation}
\label{intro-int-1-infini}
\int_1^{+\infty} e^{sy}\nu(dy)<+\infty,\qquad {\rm for\;some\;}\quad s>0\eqp{.}
\end{equation}
a ) Introducing adapted functional Banach spaces, we establish  (cf. Theorem~\ref{thm:proprietes-de-lambda} and Proposition~\ref{prop:unicite-solution}) that  $F(\theta,\mu,\rho,.)$ is the unique solution of (\ref{equa-fonctionnelle}). Moreover if $\theta>0$ or $\theta=0$ and $\mathbb E(X_1)<0$, then $F(\theta,\mu,\rho,.)$ has a sub-exponential decay~:
\begin{equation}
\label{intro-equivalent-de-F}
F(\theta,\mu,\rho,x)\leq Ce^{-\gamma x}, \quad \forall x\geq0 , \quad {\rm for\; some\;}\quad C>0, \gamma>0\eqp{.}
\end{equation}
The optimal value of $\gamma$ will be given in (\ref{renouv-existence-de-kappa}).

Note that if $\theta=0$ and $\mathbb E(X_1)\geq 0$, then $F(0,0,0,x)=1$, hence there is no hope to obtain a sub-exponential decay.

b ) Unfortunately the equation (\ref{equa-fonctionnelle}) does not permit to determine $F(\theta,\mu,\rho,.)$ explicitely, but allows to obtain an approximation scheme.  Suppose $\theta>0$ or $\theta=0$ and $\mathbb E(X_1)<0$.
We define by induction a sequence of functions $\alpha_n(\theta,\mu,\rho,.)$ verifying (\ref{intro-equivalent-de-F}) and strongly approximating $F(\theta,\mu,\rho,.)$~:
\begin{equation}
\lim_{n\rightarrow +\infty}\left(\sup_{x\geq0}|(F-\alpha_n)(\theta,\mu,\rho,x)|e^{\gamma x}\right)\leq \xi^n K\eqp{,}
\end{equation} 
where $K$ is a constant, and $\xi\in]0,1[$ and  depends on $\theta$ and $\gamma$.

\item It is worth pointing out that the  previous analysis is only valid if
\mbox{$\nu(\mathbb R)<\infty$}. To remove this assumption, we introduce the Laplace transform $\widehat F$ of $F$  with respect to the $x$ variable~: 
\begin{equation}
\widehat F(\theta,\mu,\rho,q):=\int_0^{+\infty}e^{-qx}F(\theta,\mu,\rho,x)dx\eqp{.}
\end{equation}
Since $|F|\leq1$, $\widehat F(\theta,\mu,\rho,q)$ is well defined for any $q\in \mathbb C$, $\Re q>0$.  

Suppose $\nu(\mathbb R)<+\infty$ and (\ref{intro-int-1-infini}). Starting with the integral equation (\ref{equa-fonctionnelle}) satisfied by $F(\theta,\mu,\rho,.)$, we prove that if moreover~:
\begin{equation}
\label{moment-expo-infini-1}
\int_{-\infty}^{-1}e^{-qy}\nu(dy)<+\infty,\qquad \forall q > 0\eqp{,}
\end{equation}
then  $\widehat F(\theta,\mu,\rho,.)$ verifies some kind of integral  equation (identity (\ref{equa-Laplace}) in  Theorem~\ref{thm:equa-Laplace}).
We observe that (\ref{equa-Laplace}) is still valid when $\nu$ is a Lévy measure satisfying (\ref{intro-int-1-infini}) and (\ref{moment-expo-infini-1}). Using an approximation scheme, it is easy to prove that $\widehat F(\theta,\mu,\rho,.)$ verifies (\ref{equa-Laplace}), under  (\ref{intro-int-1-infini}) and (\ref{moment-expo-infini-1}).
Moreover the  equation (\ref{equa-Laplace}) gives an equation satisfied by  the factors of the Wiener-Hopf decomposition of $\displaystyle \frac{\theta}{\theta+\varphi(-q)}$ (see for detail Remark~\ref{rem:Wiener-Hopf}, 7.).\\
In the particuler case of  the support of $\nu$ is included in $[0,+\infty[$ (i.e. $(\X)$ has only positive jumps) then $\widehat F(\theta,\mu,\rho,.)$ is explicit.
\item In section~\ref{sec:new-funct-equa} we draw a first important consequence of Theorem~\ref{thm:equa-Laplace}. We show that we can go back to $F(\theta,\mu,\rho,.)$. A simple modification in (\ref{equa-Laplace})  allows to prove that $F(\theta,\mu,\rho,.)$ verifies an integro-differential equation (cf. Theorem~\ref{thm:seconde-equa-fonc}) which is new and different from the equation verified by $F(\theta,\mu,\rho,.)$ when $\nu(\mathbb R)<+\infty$.
\item From equation (\ref{equa-Laplace}), we deduce in \cite{roynette-vallois-volpi-2003} two main consequences~:
\begin{itemize}
\item[a)]If $\nu$ has finite exponential moments, $F(\theta,\mu,\rho,.)$ has the following expansion~: 
\begin{align*}
\label{intro-gene:dev-asymp}
 F(\theta,&\mu,\rho,x) =
 C_0(\theta,\mu,\rho) e^{- \gamma_0(\theta) x} \\
& +
       \sum_{i=1}^{p}a_i \left( C_i(\theta,\mu,\rho,x) e^{-\gamma_i(\theta) x}
+\overline{C_i}(\theta,\mu,\rho,x)
       e^{-\overline{\gamma_i}(\theta)x}\right)\nonumber
+\mathrm{O}\left(e^{-Bx}\right)\eqp{,}
\end{align*}
where $C_0(\theta,\mu,\rho)$ is a positive real number,
$C_1(\theta,\mu,\rho,x), \cdots, C_p(\theta,\mu,\rho,x)$ are $x$-polynomial
functions with values in $\mathbb C$,\goodbreak 
$\left(\gamma_0(\theta),\gamma_1(\theta),
 \cdots,\gamma_p(\theta), \overline{\gamma_1}(\theta),
  \cdots, \overline{\gamma_p}(\theta)\right)$ are zeros of
   $\varphi-\theta$ (where $\varphi$ is the function defined by (\ref{def-phi}))
   and $\displaystyle a_i=\frac{1}{2}$ (resp.\ $1$)
   if  $\gamma_i(\theta)\in \mathbb R$  (otherwise).\\
 This result is an extension of the one  of J.~Bertoin and  R.A.~Doney \cite{BertoinDoney94}.
\item[b)] The asymptotic behaviour of the law of the triplet $(T_x,K_x, L_x)$ when $x\rightarrow +\infty$.
\end{itemize}
\item  Let $(\X)$ be a Lévy process. 
It is well known that there exists a family of probability measures  $\left(\mathbb P^{(\lambda)}\:,\;0\leq\lambda\leq \gamma\right)$ such that,   under $\mathbb P^{(\lambda)}$, $(\X)$ is still a Lévy process and~:
\begin{equation}
\mathbb P^{(\lambda)}(X_t\in dx)=e^{\lambda x}e^{-t\varphi(-\lambda)}\mathbb P(X_t\in dx).
\end{equation}
Consequently $\displaystyle \varphi^{(\lambda)}(q)=\varphi(q-\lambda)-\varphi(-\lambda)$, where $\varphi^{\lambda}$ is associated with $(\X)$ under $\mathbb P^{(\lambda)}$. Suppose $\nu$ verifies the assumption of Theorem~\ref{thm:equa-Laplace}, then there exists $\lambda$ such that $\varphi(-\lambda)=\theta$ and $\varphi'(0)\varphi'(-\lambda)<0$. Since $\mathbb E(X_1)=-\varphi'(0)$, and $\mathbb E^{(\lambda)}(X_1)=-\varphi'^{(\lambda)}(0)=-\varphi'(-\lambda)$, then $\mathbb E(X_1)\mathbb E^{(\lambda)}(X_1)<0$. This trick allows to only consider the case $\mathbb E(X_1)>0$ (or  \mbox{$\mathbb E(X_1)<0$}), and then simplify the proofs of Theorems  2.4,  4.1 of \cite{roynette-vallois-volpi-2003}.
\item There is a hudge litterature concerning the so-called ruin problem. A good reference for the reader interested in this topic is the book written by   T.~Rolski, H. Schmidli and J. Teugels  \cite{RolskiSchmidli99}. Historically, the first model (called classical model or the Cramér-Lundberg model) was initiated by F. Lundberg \cite{Lundberg03} and  H. Cramér \cite{Cramer30}, \cite{Cramer55}. It corresponds to the case~:\; 
$\displaystyle  X_t=-ct+J_t$, ($\J$) being a compound  Poisson process. There is no Brownian component (i.e. $\sigma=0$). The authors proved that the Laplace transform $\widehat F(0,0,0,.)$ of the ruin probability verifies a relation, and computed explicitely $\mathbb P(T_x<+\infty)$ when the jumps are exponentially distributed. Among the authors working with the classical model we  we may mention Gerber \cite{Gerber73}), F. Delbaen, J. Haezendonck \cite{DelbaenHaezendonck85}), A. Dassios, P. Embrecht \cite{DassiosEmbrechts89}, G.C. Taylor \cite{Taylor76}) and  W.  Feller \cite{Feller71}).\\
The perturbed model was introduced by H.U. Gerber \cite{Gerber70} and corresponds to our underlying process $(\X)$ with $\sigma>0$.

\noindent In some specific cases, the ruin probability, the law of $T_x$ or the distribution of the overshoot have been determined, more or less explicitely, see for instance \cite{DufresneGerber90}, \cite{KouWang01}, \cite{DozziVallois97}, \cite{HuzakPerman03}   \cite{DufresneGerberShiu91}, \cite{DicksonWaters92}, \cite{Norberg}, \cite{GerberShiu97}, \cite{NguyenYor01}. 

\end{enumerate}

\sectio{A functional relation satisfied by $F$, when  \mbox{$\nu(\mathbb R)<+\infty$}}
\label{sec:equa-fonctionnelle}


\subsection{Functional equation satisfied by $F$}


We keep the notations given in the Introduction. In this section it is assumed that~:
\begin{equation}
\label{ppc}
\lambda:=\nu(\mathbb{R})=\int_{-\infty}^{+\infty}\nu (dy)<+\infty  \eqp{.}
\end{equation}
Then ($\J$) is a compound Poisson process. Hence it admits a
first jump time  $\tau_1$, exponentially distributed with parameter
$\nu (\mathbb{R})$ and the process \mbox{$(X_{t+\tau_1}-X_{\tau_1}\,;\;t\geq0)$} is again a Lévy process distributed as ($\X$). This property is the key of our approach that we briefly describe. We distinguish three cases~:

\begin{itemize}
\item \mbox{$\displaystyle T_x:=\inf\;{\{t\geq 0\,;\; B_{t }-c_0t>x
      \}}<\tau_1$}\quad if\quad \mbox{$\displaystyle \sup_{0\leq t\leq
      \tau_1}(B_t-c_0t)>x$},
\item \mbox{$\displaystyle T_x=\tau_1$} \quad if \quad  \mbox{$\displaystyle \sup_{0\leq t\leq \tau_1}(B_t-c_0t)<x$} \quad  and \quad  \mbox{$\displaystyle J_{\tau_1}+B_{\tau_1}-c_0\tau_1 >x$},
\item \mbox{$\displaystyle T_x>\tau_1$} otherwise. However, conditionally to $\{T_x>\tau_1\}$,   $T_x-\tau_1$ is distributed as \mbox{$\displaystyle \widehat {T}_{x-X_{\tau_1}}$} where $(\widehat T_x\,;\;x>0)$ is an independent copy of $(T_x\,;\;x>0)$, independent of ($\X$).
This renewal part gives rise to the integral kernel $\Lambda_\theta$ defined by (\ref{def-de-lambda}) below.
\end{itemize}

\noindent This leads us to decompose $F(\theta,\mu,\rho,.)$ defined in (\ref{def-de-F-theta-mu-rho}), as follows~:
\begin{align}
\label{decomposition-de-F}
F(\theta,\mu,\rho,x)= & \mathbb {E}\left(e^{-\theta T_x-\mu K_x-\rho L_x}\un_{\{T_x<\tau_1\}}\right) +\mathbb {E} \left(e^{-\theta T_x-\mu K_x-\rho L_x} \un_{\{T_x=\tau_1\}}\right)\nonumber\\
                 &+ \mathbb {E} \left(e^{-\theta T_x-\mu K_x-\rho L_x}\un_{\{\tau_1<T_x<+\infty\}}\right)
\eqp{.} 
\end{align}

Finally the main result of this subsection is the following.

\begin{thm}
\label{thm:equa-fonctionnelle}
\quad Assume $\lambda=\nu(\mathbb R)<+\infty$.  For any $\theta \geq
0$, $ \mu \geq 0$ and $\rho \geq 0$, the function
$F(\theta,\mu,\rho,.)$ is solution of the following integral
equation~:
\begin {equation}
\label{equa-fonctionnelle}
G(x) =
F_0(\theta,\mu,\rho,x)+F_1(\theta,\mu,\rho,x)+\Lambda_{\theta} G(x) \qquad \forall x\ge 0
\end{equation}
where
\begin{eqnarray}
\label{alpha-theta}
\alpha_{\theta}&=&\sqrt{c_0^2+2(\lambda+\theta)}\eqp{,}\\[1ex]
\label{proba1}
F_0(\theta,\mu,\rho,x)&=& e^{-(c_0+\alpha_{\theta})x}\eqp{,} \\[1ex]
\label{proba2}
F_1(\theta,\mu,\rho,x)&=& \frac{e^{-(c_0+\alpha_{\theta})x}}{\alpha_{\theta}(\mu-\rho+c_0+\alpha_{\theta})}
           \int_{[0,x]} \left( e^{(-\rho+c_0+\alpha_{\theta})y}-e^{-\mu y} \right)\, \nu(dy)\nonumber\\[1ex] 
       && + \frac{e^{-\rho x}}{\alpha_{\theta}(\mu-\rho+c_0-\alpha_{\theta})} \; 
           \int_{]x,+\infty[}
         \left(e^{-(\rho+\alpha_{\theta}-c_0)(y-x)}-e^{-\mu (y-x)}\right) \, \nu(dy) 
                \nonumber\\[1ex]
       && + \frac{e^{(\mu-\rho) x}-e^{-(c_0+\alpha_{\theta})x}}
                    {\alpha_{\theta}(\mu -\rho+c_0+\alpha_{\theta})} \; 
           \int_{]x,+\infty[}  e^{-\mu y} \, \nu(dy)  \nonumber\\[1ex] 
       && -\frac{e^{-(c_0+\alpha_{\theta})x}}
            {\alpha_{\theta}(\mu-\rho+c_0-\alpha_{\theta})} \;
            \int_0^{+\infty}
           \left(e^{-(\rho+\alpha_{\theta}-c_0)y}-e^{-\mu y} \right)\, \nu(dy)
             \eqp{,}
\end{eqnarray}

and   $\Lambda_{\theta}$ is the operator~:
\begin{equation}
\label{def-de-lambda}
\Lambda_{\theta} G(x) 
       =\frac{1}{\alpha_{\theta}} 
           \int_{-\infty}^{+\infty} \nu(dy)\,
           \int_{-\infty}^{(x-y)\wedge x} e^{-c_0a}
             \left( e^{-\alpha_{\theta}|a|} - e^{-(2x-a)\alpha_{\theta}} \right) G(x-a-y) da \eqp{.}
\end{equation}
\end {thm}
\noindent {\bf Proof of Theorem~\ref{thm:equa-fonctionnelle}}
\\
We compute the two first terms in (\ref{decomposition-de-F}) in Lemmas~\ref{lem:proba1}, \ref{lem:proba2} and the last one in Lemma \ref{lem:proba3}.

\begin{lem}
\label{lem:proba1}
\qquad Let $\alpha_\theta$  be the real number, defined by  (\ref{alpha-theta}), then~:
\begin{equation}
\label{equa:proba1}
\mathbb {E}\left(e^{-\theta T_x -\mu K_x-\rho L_x} \un_{\{ T_x<\tau_1\}}\right)
  =e^{-(c_0+\alpha_{\theta})x} \eqp{.} 
\end{equation}
\end{lem}
\noindent{\bf Proof of Lemma~\ref{lem:proba1}}
\\
Let ($\Btild$) be the Brownian motion with drift $-c_0$~:
\begin{equation}
\label{def-Btild}
\widetilde B_t=B_t-c_0t \qquad \forall t\geq 0 \eqp{.} 
\end{equation}
We set $\widetilde T_x:= \inf{\{t\geq 0; \;\widetilde B_t > x \}}$, $x\geq0$.

\noindent Then \mbox{$\displaystyle \{T_x<\tau_1\}=\{\widetilde T_x<\tau_1\}$} and on \mbox{$\displaystyle \{T_x<\tau_1\}$}, we have 
\mbox{$\displaystyle K_x=L_x=0$}.

\noindent Since $\tau_1$ is exponentially distributed with parameter $\lambda$ and independent of $\widetilde T_x$, we have~:
\begin{equation}
\mathbb {E}\left(e^{-\theta T_x-\mu K_x-\rho L_x} \un_{\{T_x<\tau_1\}} \right)= \mathbb{E} \left(e^{-(\lambda + \theta) \widetilde T_x}\right)
\eqp{.}
\end{equation}
By (\cite{Karatszas91}, exercise $5.10$ page $197$) we can conclude that (\ref{equa:proba1}) holds.\finpreuve
\begin{lem}
\label{lem:proba2}

We have~:
\begin{align}
\label{equa:proba2}
&\mathbb {E}\left(e^{-\theta T_x -\mu K_x-\rho L_x} \un_{\{ T_x=\tau_1\}}\right)
  = \frac{e^{-(c_0+\alpha_{\theta})x}}{\alpha_{\theta}(\mu-\rho+c_0+\alpha_{\theta})}
        \int_{[0,x[}\!\!\! \left( e^{(-\rho+c_0+\alpha_{\theta})y}-e^{-\mu y} \right) \nu(dy)
             \nonumber\\[1ex] 
       &\hphantom{======}+ \frac{e^{-\rho x}}{\alpha_{\theta}(\mu-\rho+c_0-\alpha_{\theta})} \; 
           \int_{[x,+\infty[}
         \left(e^{-(\rho+\alpha_{\theta}-c_0)(y-x)}-e^{-\mu (y-x)}\right) \, \nu(dy) 
                \nonumber\\[1ex]
       &\hphantom{======}+ \frac{e^{(\mu-\rho) x}-e^{-(c_0+\alpha_{\theta})x}}
                    {\alpha_{\theta}(\mu-\rho +c_0+\alpha_{\theta})} \; 
           \int_{[x,+\infty[}  e^{-\mu y} \, \nu(dy)  \nonumber\\[1ex] 
       &\hphantom{======}-\frac{e^{-(c_0+\alpha_{\theta})x}}
            {\alpha_{\theta}(\mu-\rho+c_0-\alpha_{\theta})} \;
            \int_0^{+\infty}
           \left(e^{-(\rho+\alpha_{\theta}-c_0)y}-e^{-\mu y} \right)\, \nu(dy)
\eqp{.}\nonumber\\ 
\end{align}
\end{lem}

\noindent{\bf Proof of Lemma~\ref{lem:proba2}}
\\
Write $Y_1:=J_{\tau_1}$. We observe that on  $\{T_x=\tau_1\}$,  $Y_1>0$.
Morever~:
\begin{equation}
\{T_x=\tau_1\}=\{\sup_{t \leq \tau_1}\widetilde B_t<x, \widetilde {B}_{\tau_1 }+Y_1>x\}\eqp{,}
\end{equation}
and
\begin{equation}
K_x=\widetilde B_{\tau_1}+Y_1-x\quad,\quad \mbox{$\displaystyle L_x=x-\widetilde B_{\tau_1}$} \eqp{,}
\end{equation}
where ($\Btild$) is defined by the relation (\ref{def-Btild}). Since the distribution of $Y_1$ is $\displaystyle \frac{1}{\lambda}\nu$, conditioning by $\tau_1$ and $Y_1$, we have~:
\begin{align}
\Delta:=&\mathbb{E}\left(e^{-\theta T_x-\mu K_x-\rho L_x}\un_{\{T_x=\tau_1\}}\right)\nonumber\\
=&
e^{-\rho x}
\int_0^{+\infty} dt\; e^{-(\lambda+\theta ) t} \int_0^{+\infty}\nu(dy) \;
\mathbb{E}
\left(e^{-(\mu-\rho)\widetilde B_t-\mu(y-x)} 
\un_{\{\sup_{u \leq t} \widetilde B_u<x\;;
\;\; x-y \leq \widetilde B_t\}}\right)\nonumber\\
\label{esp-tx=tau1}
\end{align}
The density function of \mbox{$\displaystyle (\sup_{u \leq t} B_u,
 B_t)$} is given by (\cite{Karatszas91} page $95$), i.e.~:
\begin{equation}
\mathbb{P}
\left( B_t \in da; \sup_{u\leq t} B_u \in db \right)= 
\frac{2(2b-a)}{\sqrt{2 \pi t^3}}\; e^{- \frac{(2b-a)^2}{2t}} 
\un_{\{ a<b ;\, b>0\}} da db  \eqp{.}
\end{equation}
Applying Girsanov's formula, we get~:
\begin{equation}
\label{girsanov}
\mathbb{P}\left( \widetilde B_t \in da; \sup_{u\leq t}\widetilde  B_u \in db \right)= 
 \frac{2(2b-a)}{\sqrt{2 \pi t^3}}\; e^{-c_0a-\frac{c_0^2}{2}t} \;
 e^{-\frac{(2b-a)^2}{2t}} \un_{\{a<b;\,b>0\}} da db  \eqp{.}
\end{equation}
Combining (\ref{girsanov}) and (\ref{esp-tx=tau1}) leads to~:
\begin{align}
\Delta &= e^{(\mu-\rho) x} \int_0^{\infty} \nu(dy)\; e^{-\mu y}
\int_{x-y}^x da\; e^{-(c_0+\mu-\rho)a}
\int_{a \vee 0}^x db\; (2b-a)\nonumber\\[1ex]
&\hspace*{40mm}
\int_0^{\infty} \frac{2}{\sqrt{2 \pi t^3}}\; 
e^{-\frac{1}{2}\left((2(\lambda +\theta) +c_0^2)t+ \frac{(2b-a)^2}{t}\right)} dt\eqp{.}
\label{horreur}
\end{align}
%
%
%
Recall the classical identities (cf. \cite{gradshteyn80} sections  8.432 6 page 959,  and 8.469 3 page 967)~:
\begin{equation}
\mathcal{K}_{\frac{1}{2}}(\delta):=\frac{1}{2}
\int_0^{+\infty}\frac{1}{\sqrt{t}}
e^{-\frac{\delta}{2}(t+\frac{1}{t})}dt=\sqrt{\frac{\pi}{2 \delta}}e^{- \delta}  \qquad \forall \delta >0 \eqp{,}
\end{equation}
and
 \begin{equation}
\label{integrale-en-t3}
\int_0^{+\infty}\frac{1}{\sqrt{t^3}}e^{-\frac{1}{2}(\beta t+\frac{\gamma}{t})} dt = \sqrt{\frac{2\pi}{\gamma}}e^{-\sqrt{\beta \gamma}} 
\qquad \forall \beta > 0, \;\forall \gamma >0 \eqp{.}
\end{equation}
obtained  by derivation and the changing  variable  \mbox{$ t \rightarrow \sqrt{\frac{\beta}{\gamma}} \; t$}. 

\noindent This allows to  first compute explicitely the integral with respect to $dt$ in (\ref{horreur})~: 
\begin{equation}
\Delta=
2e^{(\mu-\rho) x}
\int_0^{+\infty} \nu(dy)\; e^{-\mu y}
\int_{x-y}^x da\; e^{-(c_0+\mu-\rho)a}
\int_{a \vee 0}^x e^{-\alpha_{\theta}(2b-a)}\; db \qquad \qquad \qquad
\end{equation}
 In a second step we evaluate the integral with respect to $db$~:
\begin{equation}
\Delta=
\frac{e^{(\mu-\rho) x}}{\alpha_{\theta}}
\int_0^{+\infty} \nu(dy)\;  e^{-\mu y}
\int_{x-y}^x da\; e^{-(c_0+\mu-\rho)a}
\left( e^{-\alpha_{\theta}|a|}- e^{-\alpha_{\theta}(2x-a)} \right)\eqp{.}
\end{equation}
To drop $|a|$, we introduce two cases  $x-y\geq 0$ and  $x-y<0$~:
\begin{align}
\Delta=&
\frac{e^{(\mu-\rho) x}}{\alpha_{\theta}}
\int_{[0,x[} \nu(dy)\; e^{-\mu y}
\int_{x-y}^x  e^{-(c_0+\mu-\rho+\alpha_{\theta})a}\; da\nonumber\\[1ex]
&+
\frac{e^{(\mu-\rho) x}}{\alpha_{\theta}}
\int_{[x,+\infty[} \nu(dy)\; e^{-\mu y}
\left[\int_{x-y}^0  e^{-(c_0+\mu-\rho-\alpha_{\theta})a}\; da +
\int_0^x  e^{-(c_0+\mu-\rho+\alpha_{\theta})a}\; da \right]
\nonumber\\[1ex]
&-
\frac{e^{(-2\alpha_{\theta}+\mu-\rho) x}}{\alpha_{\theta}}
\int_0^{+\infty} \nu(dy)\; e^{-\mu y}
\int_{x-y}^x  e^{-(c_0+\mu-\rho-\alpha_{\theta})a}\;da \eqp{.}
\end{align}
Computing the integral with respect to $da$ we easily obtain (\ref{equa:proba2}).\finpreuve

\begin{lem}
\label{lem:proba3}
\quad In (\ref{decomposition-de-F}), the third expectation is equal to~:
\begin{align}
\label{equa:proba3}
&\mathbb{E}\left( e^{-\theta T_x-\mu K_x-\rho L_x}\; \un_{\{\tau_1<T_x<+\infty\}}\right)=
\nonumber\\[1ex]
& {\displaystyle\frac{1}{\alpha_\theta}}\!\!\!
\int_{-\infty}^{+\infty}\!\!\!\!\!\! \nu(dy) 
\int_{-\infty}^{(x-y) \wedge x}\!\!\!\!\!e^{-c_0a}
\left(e^{-\alpha_\theta |a|}-e^{(2x-a)\alpha_\theta}\right)
F(\theta,\mu,\rho,x-a-y) da\eqp{.}
\end{align}
Morever~:
\begin{equation}
\label{operateur-positif}
e^{-\alpha_\theta|a|}-e^{-(2x-a)\alpha_\theta}\geq0\qquad {\rm if}\qquad a\leq (x-y)\wedge x\eqp{,}
\end{equation}
 so $\Lambda_\theta$ défined by  (\ref{def-de-lambda}) is an non-negative operator.
\end{lem}

\noindent{\bf Proof of Lemma~\ref{lem:proba3}}
\\
Formula (\ref{equa:proba3}) may be proved proceding analogously to the proof of previous Lemma.\finpreuve


\subsection{Study of  $\Lambda_\theta$}
\label{subsec:operateur-lambda}


To investigate uniqueness in (\ref{equa-fonctionnelle}), we prove that 
$\Lambda_\theta$ is a contraction in some functional Banach spaces. Recall we have in mind to prove  that $F$ has a sub-exponential decay at infinity, therefore it seems natural to introduce the Banach space~:
\begin{equation}
\label{Banach}
\mathcal{B}_\gamma:=\{ f :\mathbb R_+ \rightarrow \mathbb R \, ;
\,\sup_{ x \in \mathbb R_+}e^{\gamma \, x}| f(x) | < +\infty\ \}\quad \gamma\geq0 \eqp{.}
\end{equation}
$\mathcal{B}_\gamma$ is equipped with the norm~:
\begin{equation}
\label{norme}
\|f\|_\gamma=\sup_{x \in \mathbb R_+}e^{\gamma\,x}|f(x)|\eqp{.}
\end{equation}
The values of $\gamma$ such that  $\Lambda_\theta$ is a contraction in  $\mathcal{B}_\gamma$  are linked to the zeros of the function \mbox{$\displaystyle \varphi_\theta$}~:
\begin{equation}
\label{def-phi-theta}
\varphi_\theta(q):=\varphi(q)-\theta \eqp{.}
\end{equation}
Before stating our main result, we fix some notations. Let~:
\begin{align}
\label{def-de-r}
&r_\nu := \sup \left\{s \geq 0 ;\; \int_1^{+\infty} e^{sy}
\nu (dy )\;<\;+\infty \right\}\eqp{,}
\end{align} 
with the convention $\sup \emptyset =0$.\\

 $\widehat {\nu}$ (resp.\ $\widehat{\nu^+} $) denotes the Laplace transform of $\nu$ (resp.\ $\displaystyle \nu_{|_{[0,+\infty[}}$), if they exist~:
\begin{align}
\widehat{\nu}(q):= &\int_{- \infty}^{+ \infty}e^{-qy} \nu(dy) 
\\
\label{def-de-trans-nu+}
\widehat{\nu^+}(q):=& \int_0^{+ \infty}e^{-qy} \nu(dy) \eqp{.}
\end{align}

\begin {thm}
\label {thm:proprietes-de-lambda}
\quad Suppose that $\nu(\mathbb R)<+\infty$.
\begin {itemize}
\item [(i)] For any $ \theta \geq 0$, the operator $\Lambda_\theta$
defined by (\ref{def-de-lambda}) is a linear and non-negative operator, with norm equals to  $\displaystyle \frac{\lambda}{\lambda+\theta}$ in $L^{\infty}(\mathbb{R_+})$.
\item [(ii)]Assume  $r_\nu >0$. Let $ \gamma \in [0, r_\nu[$
and  $\theta > 0$ or $\theta=0$ and $\mathbb E(X_1)<0$. Then~:
\begin{itemize}
\item[a)]
$\Lambda_\theta$ is a bounded operator from 
$\mathcal{B}_\gamma$ to $\mathcal{B}_\gamma$. More precisely~:
\begin{equation}
\label{Lambdaborne}
\| \Lambda_\theta f \|_\gamma \le c_{\theta,\gamma}\|f\|_\gamma \quad
\forall f \in \mathcal B_\gamma \eqp{,}
\end{equation}
with
\begin{equation}
\label{c-theta-gamma} c_{\theta,\gamma}
=\frac {\widehat{\nu}(-\gamma)} {\widehat{\nu}(-\gamma)-
\varphi(-\gamma)+\theta}\eqp{.}
\end{equation}
\item[b)] There exists $\gamma \in ]0,r_\nu[$ such that $\varphi(-\gamma)<\theta$. Therefore $\Lambda_\theta$ is a contraction in $\mathcal{B}_\gamma$ since~:
\begin{equation}
\label{majoration-cgamma}
0 < c_{\theta,\gamma} < 1\eqp{.}
\end{equation}
\end{itemize}
\end{itemize}
\end{thm}
\begin{rem}
\label{rem:intervalle-important}
 {\rm
\begin{enumerate}
\item It is clear  that  $r_\nu >0$ is equivalent to (\ref{intro-int-1-infini}) and if \mbox{$\displaystyle r_\nu \in]0;+\infty]$}, then $\varphi$ given formally by (\ref{def-phi}),  is actually well-defined  on  $]-r_\nu,0]$.
\item If $r_{\nu} >0$, it is easy to check 
(cf.  Annex~\ref{A1-fonction-phi}) that   \mbox{$\displaystyle \{q \in\; ]-r_\nu,0]\; / \;\varphi(q)<\theta\}$} is non empty if   $\theta>0$ or if  $\theta=0$ and \mbox{$\mathbb E(X_1)<0$}.
Observe that $X_1$ has a finite expectation if~:
\begin{equation}
\int_{\mathbb R}|y|\nu(dy)<+\infty\eqp{.}
\end{equation}
In this case~:
\begin{equation}
\mathbb E(X_1)=-c_0+\int_{\mathbb R}y\nu(dy)=-c+\int_{\mathbb R}y\un_{\{|y|>1\}}\nu(dy)\eqp{,}
\end{equation}
where $c$ comes from Lévy-Khintchine formula (\ref{def-phi}).
\item The assumption (\ref{intro-int-1-infini}) means that the positive jumps are not too big. It corresponds to the intuition, since more the positive jumps are small, more time is needed to reach a positive level $x$. Hence more $F(\theta,\mu,\rho,.)$ decreases.
\end{enumerate}
}
\end{rem}

\noindent {\bf Proof of Theorem~\ref{thm:proprietes-de-lambda}} 

\noindent $(i)$ Relation (\ref{operateur-positif}) implies that $\Lambda_\theta$ is a non-negative operator.\\
It is easy to check that the function $\ell$~: 
\begin{equation} 
\ell( x )= \un_{\{a+y \leq x\,;\; a\leq x\}}e^{-c_0a}
\left(e^{-\alpha_\theta |a|}-e^{-(2x-a)\alpha_\theta}\right)
\end{equation} 
is increasing, then~:
\begin{equation}
\label{majo-ell}
\ell(x) <\ell(+\infty)= e^{-c_0a} e^{-\alpha_\theta
|a|}\qquad \forall x \in \mathbb{R}\eqp{.}
\end{equation}
A straightforward calculation shows that  \mbox{$\displaystyle |\Lambda_\theta h(x)|\leq 
\frac{\lambda}{\lambda+\theta}\;\|h\|_{\infty}$}, for any \mbox{$x\geq0$}.

\noindent Taking \mbox{$\displaystyle h:x\rightarrow 1$}, we have \mbox{$\displaystyle \|\Lambda_\theta h\|_{\infty}=\frac{\lambda}{\lambda+\theta}$}, then
\mbox{$\displaystyle \||\Lambda_\theta|\|_{L^{\infty}(\mathbb {R_+})}=
\frac{\lambda}{\lambda+\theta}$}.

\noindent $(ii)$ Let $f$ be an element of  $\displaystyle \mathcal {B}_\gamma$, then  \mbox{ $\displaystyle  |f(x)|\leq \|f\|_\gamma\; e^{-\gamma x}$}, $\forall x \geq 0$. Consequently~:
\begin{equation}
|\Lambda_\theta f(x)| \leq
\frac{1}{\alpha_\theta}\; \|f\|_{\gamma} \; e^{-\gamma x}
\int_{-\infty}^{+\infty}\!\!\!\!\!\! \nu(dy)
\int_{-\infty}^{(x-y)\wedge x}\!\!\!\!\!\!\!\!\! \! e^{-c_0a}\left(e^{-\alpha_\theta |a|}-e^{-(2x-a)\alpha_\theta}\right)e^{\gamma(a+y)}da\eqp{,} 
\end{equation}
for any  $\gamma \in \;[0,r_\nu[$.

\noindent Making use of (\ref{majo-ell}), we get~:
\begin{equation}
|\Lambda_\theta f(x)| 
\leq 
\frac{1}{\alpha_\theta}\; \|f\|_{\gamma} \; e^{-\gamma x}\!
\int_{-\infty}^{+\infty}\!\!\! \nu(dy)e^{\gamma y}
\left[
\int_{-\infty}^0\!\! e^{-(c_0-\alpha_\theta -\gamma)a} da +
\int_0^{+\infty}\!\!\! e^{-(c_0+\alpha_\theta -\gamma)a} da \right]\nonumber\\[1ex]
\end{equation}
Computing the integral with respect to $da$, yields directly to (\ref{Lambdaborne}). \finpreuve

\begin{prop}
\label{prop:unicite-solution} \quad 
Assume $\nu(\mathbb R)<+\infty$, $r_\nu >0$,  $\mu\geq 0$,  
$\theta> 0$ or   $\theta=0$  if  $\mathbb E(X_1)<0$. Let  $\gamma$ be in  $[0,r_\nu[$, such that $\varphi(-\gamma)< \theta$. Then the  function $ F(\theta,\mu,\rho,.)$ belongs to  $\mathcal B_\gamma$ and
the equation  (\ref{equa-fonctionnelle}) has an  unique
solution in $\mathcal B_\gamma$.
\end{prop}

\noindent To prove  Proposition~\ref{prop:unicite-solution}, we need the following preliminary.

\begin{lem}
\label{lem:limite-lambda-n} \quad 
Suppose $\theta> 0$ or   $\theta=0$  if  $\mathbb E(X_1)<0$, then for any $x>0$,
\begin{equation}
\label{limite-lambda-n}
\lim_{n \rightarrow +\infty} \Lambda_\theta^{n} F(\theta,\mu,\rho,.)(x)=0\eqp{.}
\end{equation}
\end{lem}

\paragraph{ Proof of Lemma~\ref{lem:limite-lambda-n}}\mbox{}

\noindent $1)$  Suppose $\theta>0$. Since $F$ is bounded by $1$, and the norm of $\Lambda_\theta$ is $\displaystyle \frac{\lambda}{\lambda+\theta}$ (cf. Theorem~\ref{thm:proprietes-de-lambda})~:
\mbox{$\displaystyle \|\Lambda_\theta^n F(\theta,\mu,\rho,.)\|_{\infty}\leq\left(\frac{\lambda}{\lambda+\theta}\right)^n$}. This proves (\ref{limite-lambda-n}).\\

\noindent  $2)$  We now turn to the case $\theta=0$ and $\mathbb E(X_1)<0$.
Iterating the functional equation (\ref{equa-fonctionnelle}), we obtain~:
\begin{equation}
\label{egalite-serie-suite}
F(\theta,\mu,\rho,x)=\sum_{p=0}^{n-1}\Lambda_\theta^p\left[F_0+F_1)(\theta,\mu,\rho,.)\right](x)+\Lambda_\theta^n F(\theta,\mu,\rho,.)(x)\eqp{.}\\[2ex]
\end{equation}
 The norm of $\Lambda_\theta$ in $\mathcal B_\gamma$  is strickly less than $1$, then the series in (\ref{egalite-serie-suite}) converges. Consequently the remaining term $\Lambda_\theta^n F(\theta,\mu,\rho,.)(x)$ converges in $\mathcal B_\gamma$ to some function  $G(\theta,\mu,\rho,x)$. It is easy to check the following~:
\begin{itemize}
\item [a)] $\displaystyle G(0,\mu,\rho,.)$ is a bounded and non negative function \eqp{,}
\item[b)] $\displaystyle G(0,\mu,\rho,.)$ is a continuous function on $[0,+\infty[$ \eqp{,}
\item[c)] \mbox{$\displaystyle \lim_{x \rightarrow +\infty}
G(0,\mu,\rho,x)=0$} \eqp{,}
\item[$d)$] $\displaystyle \Lambda_0 G(0,\mu,\rho,.)=G(0,\mu,\rho,.)$\eqp{.}
\end{itemize}
Hence by point $(i)$ of Theorem~\ref{thm:proprietes-de-lambda}~:
\begin{equation}
\label{der-inegal}
G(0,\mu,\rho,x)=\Lambda_0 G(0,\mu,\rho,.)(x)\leq \|G(0,\mu,\rho,.)\|_\infty, \quad x\geq0\eqp{.}
\end{equation}
As (\ref{majo-ell}) is a strict inequality then  (\ref{der-inegal}) is strict if \mbox{$\displaystyle\|G(0,\mu,\rho,.)\|_\infty \neq0$}.

 \noindent According to $b)$ and $c)$, there exists $x_0\geq 0$ such that~:
 
\mbox{$\displaystyle G(0,\mu,\rho,x_0)=\|G(0,\mu,\rho,.)\|_\infty $}. This implies $\|G(0,\mu,\rho,.)\|_\infty =0$.\finpreuve

\paragraph{Proof of Proposition~\ref{prop:unicite-solution}}\mbox{}

 Using the explicit expression of $F_0$ and $F_1$ (cf. (\ref{proba1}) and (\ref{proba2})),by  a straightforward calculation, it may be concluded that \mbox{$\displaystyle F_0(\theta,\mu,\rho,.)$} and \mbox{$\displaystyle F_1(\theta,\mu,\rho,.)$}  belong to \mbox{$ \displaystyle \mathcal B_\gamma$} (for a detailed proof, cf. \cite{volpi03}).\\
 By  Lemma~\ref{lem:limite-lambda-n} and (\ref{egalite-serie-suite}),
\begin{equation}
F(\theta,\mu,\rho,x)=\sum_{n=0}^{+\infty} \Lambda_\theta^n \left(F_0+F_1)(\theta,\mu,\rho,.)\right)(x)\eqp{,}
\end{equation}
Because $\displaystyle F_0+F_1 \in \mathcal B_\gamma$ and $\Lambda_\theta$ is a contraction in $ \mathcal B_\gamma$, the serie converges in  $\mathcal B_\gamma$. This directly implies the result.\finpreuve

\begin{rem}\mbox{}\\
{\rm 
\noindent$1.$  Under the conditions stated in Proposition~\ref{prop:unicite-solution}, we have actually proved that $F(\theta,\mu,\rho,x)$ can be approximated by \mbox{$\displaystyle \sum_{n=0}^{p} \Lambda_\theta^n\left(F_0+F_1)(\theta,\mu,\rho,.)\right)(x)$}. More precisely~:
\begin{equation}
\label{approximation-numerique}
\left|F(\theta,\mu,\rho,x)-\sum_{n=0}^{p} \Lambda_\theta^n \left(F_0+F_1)(\theta,\mu,\rho,.)\right)(x)\right| < c_{\theta,\gamma}^{p+1} K\; e^{-\gamma x}\eqp{,}
\end{equation}
where  \mbox{$\displaystyle K=\left\|
\sum_{n=0}^{+\infty} \Lambda_\theta^n\left(F_0+F_1)(\theta,\mu,\rho,.)\right)\right\|_\gamma <+\infty$} and  $c_{\theta,\gamma}$ is defined by (\ref{c-theta-gamma}).\\
\noindent $2.$  Let us consider the case where the support of $\nu$ is included in $]-\infty,0]$. Then $\varphi$ is well defined on $]-\infty,0]$ and $r_\nu=+\infty$. Morever $K_x=L_x=0$ and $F_1(\theta,\mu,\rho,x)=0$ for any $x\geq0$. As a result (\ref{equa-fonctionnelle}) reduces to~:
\begin{equation}
F(\theta,\mu,\rho,x)=e^{-(c_+\alpha_\theta)x}+\Lambda_\theta F(\theta,\mu,\rho,.)(x)\eqp{.}
\end{equation}
If $\theta>0$ or $\theta=0$ and $\mathbb E(X_1)<0$, we prove in Annex, Properties~\ref{prop:etude-de-phi} the existence of an unique real number $\gamma_0(\theta)$ such that~:
\begin{equation}
 -\gamma_0(\theta)<0 \qquad \rm{et}\qquad  \varphi(-\gamma_0(\theta))=\theta
\eqp{.}
\end{equation}
A direct (but fastidious !)  calculation   shows that \mbox{$\displaystyle
x\rightarrow e^{-\gamma_0(\theta)x}$} is a solution of
(\ref{equa-fonctionnelle}). For more details we refer the reader to \cite{volpi03}.\\
Hence $F(\theta,\mu,\rho,x)=F(\theta,0,0,x)=e^{-\gamma_0(\theta)x}$. 
}
\end{rem}

\sectio{The Laplace transform of $F(\theta,\mu,\rho,.)$}
\label{sec:equa-Laplace}
\subsection{The  Laplace transform expression $\widehat F(\theta,\mu,\rho,.)$}
\label{subsec:expression-trans-Laplace}

In the previous section we have proved that $F(\theta,\mu,\rho,.)$ verifies the integral equation (\ref{equa-fonctionnelle}) when $\nu$ is a probability measure. If moreover $r_\nu<+\infty$, then $F(\theta,\mu,\rho,.)$ is the unique solution of (\ref{equa-fonctionnelle}). Unfortunately we cannot define the operator $\Lambda_\theta$ if $\nu$ is not a probability measure. We would like to consider Lévy processes that do not reduce to a Brownian motion with drift plus a compound Lévy process.\\
Our approach is based on the use of the Laplace transform of $F(\theta,\mu,\rho,.)$. Since $F(\theta,\mu,\rho,.)$ is a bounded function on $[0,+\infty[$, its Laplace transform~:
\begin{equation}
\widehat {F}(\theta,\mu,\rho,q) : = \int_0^{+ \infty} e^{-qy} F(\theta,\mu,\rho,y) dy\eqp{,}
\end{equation}
is well defined for any $q$ such that $\Re(q)>0$.\\
We first suppose that $\nu$ is a finite measure. Taking the Laplace transform in (\ref{equa-fonctionnelle}), it is proved (cf. Theorem~\ref{thm:equa-Laplace}) that under some additional assumption, $\widehat F(\theta,\mu,\rho,.)$ verifies some kind of integral equation. In the calculations, cancellations occur so that in the final identity 
 and  $\nu(\mathbb R)<+\infty$ may be removed.\\
Before stating the main result of this sub-section (i.e. Theorem~\ref{thm:equa-Laplace}), we introduce~:
\begin{equation}
\label{def-de-D0}
D_0 = \{q \in \mathbb {C} \;; \; \Re q > 0 \} \eqp{.}
\end{equation} 
We suppose~:
\begin{equation}
\label{int-sur-R-}
\int_{-\infty}^{-1}e^{-qy}\nu(dy)\;<\; \infty \qquad \forall q>0.
\end{equation}
Let $R$ the operator~: 
\begin{equation}
\label{def-de-R}
Rh(q) : = \int_{-\infty}^0  \nu(dy) \; \int_0^{-y}\left( e^{-q(b+y)}-1\right) h(b) db  \eqp{,}
\end{equation}
where $q \in D_0$ and $h \in L^{\infty} (\mathbb {R_+})$.\\
\noindent Property (\ref{int-sur-R-}) implies  that  $Rh(q)$ exists.\\
We suppose moreover  that $r_\nu\in]0;+\infty]$,  with  $r_\nu$  defined in (\ref{def-de-r}), i.~e., there exist some $s>0$  such that $\displaystyle \int_1^{+\infty}e^{sy}\nu(dy)<+\infty$. Recall that from Remark~\ref{rem:intervalle-important}, $\varphi$ is well defined on $]-r_\nu,0]$.\\ 
We only consider the cases given by Figures \ref{fig:profit-net}$a$, \ref{fig:non-profit-net}$a$ and \ref{fig:egalite} (cf Annex). As a result, there exists  $\kappa >0$ such that, for all $\theta \in[0,\kappa]$~:
\begin{equation}
\label{renouv-existence-de-kappa}
\exists\; -\gamma_0(\theta)\in ]-r_\nu,0] \quad 
{\rm such\; that}\quad 
\varphi(-\gamma_0(\theta))=\theta\eqp{.} 
\end{equation}
More precisely~:
\begin{alignat}{2}
\label{renouv-existence-de-kappa-point-i}
&(i) \qquad  {\rm if}\quad \theta > 0,&& \qquad -\gamma_0(\theta)< 0\eqp{,} \\[1ex]
&(ii)
\label{renouv-existence-de-kappa-point-ii}
  \qquad {\rm if} \quad \theta=0 \quad {\rm et}\quad \mathbb E(X_1)<0, &&\qquad -\gamma_0(0)<0\eqp{,}\\[1ex]
&(iii)
\label{renouv-existence-de-kappa-point-iii}
 \qquad {\rm if}  \quad \theta=0 \quad {\rm et}\quad  \mathbb E(X_1) \geq 0,&&\qquad -\gamma_0(0)=0
\eqp{.}
\end{alignat}
These hypotheses are in force in the whole paragraphs 2 and 3.\\
\goodbreak  
\begin{thm}\mbox{}
\label{thm:equa-Laplace}
 We suppose $r_\nu<+\infty$, and (\ref{int-sur-R-}), (\ref{renouv-existence-de-kappa}) hold.
\begin{enumerate}
\item There exist $\gamma^*_0(\theta)$ such that
\begin{itemize}
\item If  $\theta>0$ or $\theta=0$ and $\mathbb E(X_1)>0$, 
$\gamma^*_0(\theta)$ is the unique positive real number such that~:
\begin{equation}
\label{theta-npn-gamma0*}
\gamma^*_0(\theta)>0 \qquad {\rm et}\qquad  
\varphi(\gamma^*_0(\theta))=\theta \eqp{.}
\end{equation}
\item If  $\theta=0$ and  $ \mathbb E(X_1)\leq0$ then \mbox{$\displaystyle \gamma_0^*(0)=0$}.
\end{itemize}
\item Let $\theta,\mu,\rho \geq0$, $q\in D_0$. We have~: 
\begin{align}
\widehat {F}(\theta,\mu,\rho,q) = \frac{1}{\varphi (q)-\theta}& \left (
\frac{q-\gamma^*_0(\theta)} {2} +  
\int_0^{+ \infty} \left[\frac{e^{-(q+\rho)y}-e^{-\mu y}}{q+\rho-\mu}\right. \right.
\nonumber\\[1ex]
&\left.
- \frac{e^{-(\gamma^*_0(\theta)+\rho)y}-e^{-\mu y}}{\gamma^*_0(\theta)+\rho-\mu}\right]\;  \nu(dy)\nonumber\\[1ex]
& \left.+ RF(\theta,\mu,\rho,.)(q)- RF(\theta,\mu,\rho,.)(\gamma^*_0(\theta)) \vphantom{\int}\right )
\label{equa-Laplace}  
\end{align}
where $\theta, \mu, \rho \geq0$, $q\in D_0$.
\end{enumerate} 
\end{thm}

\begin{rem}\mbox{}
\label{rem:Wiener-Hopf}
{\rm
\begin{enumerate}
\item  In (\ref{equa-Laplace}),  \mbox{$\displaystyle \frac{e^a-e^b}{a-b}$} stands for $e^a$  when $a=b$.
\item Assumptions (\ref{int-sur-R-}) and (\ref{renouv-existence-de-kappa}) are needed to obtain the existence of $\gamma^*_0(\theta)$, for any $\theta\geq0$ (cf. Annex~\ref{A1-fonction-phi}, subsection~\ref{phi-theta-etude}).
\item
  The function  $\widehat F(\theta,\mu,\rho,.)$ being defined on $ D_0$, then $q=\mu$ and  $q=\gamma^*_0(\theta)$ are false singularities of the right-hand side of (\ref{equa-Laplace}) .
\item  
If  $\nu(]-\infty,0])=0$, then  $ RF(\theta,\mu,\rho,.)$ cancels, and $ \widehat{F}(\theta,\mu,\rho,q)$ is given by the following explicit formula~:
\begin{align}
\label{trans-nu+}
\widehat {F}(\theta,\mu,\rho,q) = &\frac{1}{\varphi (q)-\theta} \left (
\frac{q{-}\gamma^*_0(\theta)} {2} +   
\int_0^{+ \infty}\left[\frac{e^{-(q+\rho)y}{-}e^{-\mu y}}{q+\rho-\mu} \right.\right.
\nonumber\\[1ex]
&
\left.\left. - \;\frac{e^{-(\gamma^*_0(\theta)+\rho)y}-e^{-\mu y}}{\gamma^*_0(\theta)+\rho-\mu} \right]   \nu(dy)
 \right )   
\end{align}
\item If $\nu(]0,+\infty[)=0$, (\ref{equa-Laplace}) is equivalent to~: 
\begin{align}
\label{trans-nu-}
\widehat {F}(\theta,\mu,\rho,q) = \frac{1}{\varphi (q){-}\theta} \left (
\frac{q{-}\gamma^*_0(\theta)} {2}
{+}RF(\theta,\mu,\rho,.)(q){-}RF(\theta,\mu,\rho,.)(\gamma^*_0(\theta))\vphantom{\int}\!\! \right)\end{align}
\item Let us detailed the case $\theta=\mu=\rho=0$ (i.e. $F(0,0,0,x)$ is the ruin probability). If $\mathbb E(X_1)\geq0$, it is easy  to check that $f:x\rightarrow 1$ verifies (\ref{equa-Laplace}).

Let us concentrate on the more interesting case~: $\mathbb E(X_1)<0$.

\noindent Relation (\ref{equa-Laplace}) becomes~:
\begin{equation}
\label{trans-laplace-ruine-profit-net}
\widehat F(0,0,0,q)=
 \frac{1}{\varphi(q)}
\left(\frac {q}{2}+\frac{1}{q} \int_0^{+\infty}\!\!\!\!\!(e^{-qy}-1+qy)\nu (dy)
+
RF(0,0,0,.)(q) \!\!\right)
\end{equation}
Suppose morever that $\nu(]-\infty,0[)=0$, then (\ref{trans-laplace-ruine-profit-net}) reduces to~:
\begin{equation}
\label{F00q+}
\widehat F(0,0,0,q)=
 \frac{1}{q}+\frac {\mathbb E (X_1)}{\varphi(q)}\eqp{.}
\end{equation}
It can be proved (see \cite{volpi03}, for details) that (\ref{F00q+}) generalizes  identity (2.9) of  \cite{DufresneGerber90}.

\noindent However if  $\nu(]0,+\infty[)=0$, then \mbox{$ \displaystyle F(0,0,0,x)=e^{-\gamma_0(0) x}$}.

\item Recall the Wiener-Hopf decomposition (cf. \cite{Bertoin96}, page 165)~: for any  $\theta>0$, we have~:
\begin{equation}
\label{Wiener-Hopf+}
\frac{\theta}{\theta+\varphi(-q)}=\psi_{\theta}^+(q)\psi_{\theta}^-(q)\eqp{,}
\end{equation}
where 
\begin{equation}
\psi_{\theta}^+(q)=\mathbb E\left(e^{iq S_{\tau_\theta}}\right)\quad,\qquad 
\psi_{\theta}^-(q)=\mathbb E\left(e^{iq (S_{\tau_\theta}-X_{\tau_\theta})}\right)\eqp{,}
\end{equation}
and $\tau_{\theta}$ is an   exponential  r.~v.\ with parameter $\theta$, independent from  process ($\X$) and  $S_t:=\sup_{s\leq t}X_s$.\\
Since~:
\begin{equation}
\mathbb P(S_{\tau_\theta}>a)=\mathbb P(T_a<\tau_\theta)=\mathbb E (e^{-\theta T_a})=F(\theta,0,0,a)\eqp{,}
\end{equation}
it is easy to deduce the following identity~:
\begin{equation}
\label{psi+}
\psi_{\theta}^+(q)=1+i q \widehat F(\theta,0,0,iq)\eqp{.}
\end{equation}
Equation  (\ref{equa-Laplace}) implies that the Wiener-Hopf factor $\psi_{\theta}^+$  verifies a functional equation. In particular, if $\nu(]-\infty,0])=0$, combining  equations  (\ref{trans-nu+}) and (\ref{psi+}) an  explicit form of $\psi_{\theta}^+(q)$ may be obtained. Due to (\ref{Wiener-Hopf+}), $\psi_{\theta}^-(q)$ is also explicit.
\end{enumerate}
}
\end{rem}

\paragraph{Proof of Theorem~\ref{thm:equa-Laplace}}\mbox{} \\

For simplicity, we prove (\ref{equa-Laplace}) in the particular case $\rho=0$, and we write $F(\theta,\mu,x) $ instead of $F(\theta,\mu,\rho,x)$.
The proof will be divided into two steps. We first prove (\ref{equa-Laplace}) when $\nu$ satisfied the assumptions given in Theorem~\ref{thm:equa-Laplace} and  \mbox{$\nu(\mathbb R)<+\infty$}. In a second step, we approximate $\nu$ by a sequence of finite measures $(\nu_n)$  and we take the limit in (\ref{equa-Laplace}).

\paragraph{  Step 1 }\mbox{} We suppose $\nu(\mathbb R)<+\infty$, $r_\nu<+\infty$ and (\ref{int-sur-R-}), (\ref{renouv-existence-de-kappa}) and (\ref{moment1-nu}) hold.

a) Taking the Laplace transform in functional equation  (\ref{equa-fonctionnelle}) leads to~:
\begin{equation} 
\label{trans-laplace-F}
\widehat{F}(\theta,\mu,q)=\widehat{F_0}(\theta,\mu,q)+\widehat{F_1}(\theta,\mu,q)+\widehat{\Lambda_\theta F}(\theta,\mu,.)(q) \qquad q\;\in \; D_0\eqp{.}
\end{equation}
Relation (\ref{proba1}) implies~:
\begin{equation} 
\label{trans-laplace-F0}
\widehat{F_0}(\theta,\mu,q)=\frac{1}{c+\alpha_\theta+q}\eqp{.}
\end{equation}
As for  $\widehat F_1(\theta,\mu,x)$, starting from (\ref{proba2}), we split the integral in four parts~:
\begin{align}
\widehat{F_1}(\theta,\mu,q)&=\int_0^{+\infty} e^{-qx}F_1(\theta,\mu,x)dx\nonumber\\[1ex]
&= I_1(\theta,\mu,q)+I_2(\theta,\mu,q)+I_3(\theta,\mu,q)+I_4(\theta,\mu,q)
\label{somme-integrales}
\end{align}
where
\begin{align}
I_1(\theta,\mu,q)&= \int_0^{+\infty} 
\frac{e^{-(q+c+\alpha_\theta)x}}{\alpha_\theta(\mu+c+\alpha_\theta)}
\left(\int_{[0,x]} \left(e^{(c+\alpha_\theta)y}-e^{-\mu y}\right) \nu(dy)\right)dx \nonumber\\[1ex]
&=\frac{1}{\alpha_\theta(q+c+\alpha_\theta)(\mu+c+\alpha_\theta)}
\int_0^{+\infty} \left(e^{-qy}-e^{-(q+\mu+c+\alpha_\theta)y}\right)\nu(dy)
\eqp{,}
\end{align}
\begin{align}
I_2(\theta,\mu,q)
&= \int_0^{+\infty} 
\frac{e^{-qx}}{\alpha_\theta(\mu+c-\alpha_\theta)}
\left(\int_{]x,+\infty[} 
\left(e^{-(\alpha_\theta-c)(y-x)}-e^{-\mu (y-x)}\right) \nu(dy)
\right)dx \nonumber\\[1ex]
&=\frac{1}{\alpha_\theta(\mu+c-\alpha_\theta)}
\int_0^{+\infty} 
\left(
\frac{e^{-(\alpha_\theta-c)y}-e^{-qy}}{q+c-\alpha_\theta}+
\frac{e^{-qy}-e^{-\mu y}}{q-\mu}
\right) \nu(dy)
\eqp{,}
\end{align}
\begin{align}
I_3(\theta,\mu,q)
&= \int_0^{+\infty} 
\frac{e^{-(q-\mu)x}-e^{-(q+c+\alpha_\theta)x}}{\alpha_\theta(\mu+c+\alpha_\theta)}
\left(\int_{]x,+\infty[} e^{-\mu y}\;\nu(dy)\right)dx
\nonumber\\[1ex] 
&=
-\frac{1}{\alpha_\theta(\mu+c+\alpha_\theta)}
\int_0^{+\infty}\!\!\!\!\!\!\!  \nu(dy)
\left(
\frac{e^{-qy}-e^{-\mu y}}{q-\mu}-
\frac{e^{-(q+\mu+c+\alpha_\theta)y}-e^{-\mu y}}{q+c+\alpha_\theta}
\right),
\end{align}
\begin{align}
I_4(\theta,\mu,q)
&= -\; \int_0^{+\infty} 
\frac{e^{-(q+c+\alpha_\theta)x}}{\alpha_\theta(\mu+c-\alpha_\theta)}\; dx 
\left(\widehat {\nu^+}(\alpha_\theta-c)-\widehat {\nu^+}(\mu)\right)
\nonumber\\[1ex] 
&=
\frac{\widehat {\nu^+}(\mu)-\widehat {\nu^+} (\alpha_\theta-c)}
{\alpha_\theta(q+c+\alpha_\theta)(\mu+c-\alpha_\theta)}
\eqp{,}\hphantom{====================}
\end{align}
with~:
\begin{equation*}
\widehat{\nu^+}(q):=\widehat \nu_{|_{]0,+\infty]}}(q)= \int_0^{+ \infty}e^{-qy} \nu(dy) \eqp{.}
\end{equation*}
Consequently~:
\begin{align}
&\widehat{F_1}(\theta,\mu,q)= - 
\frac{\widehat{\nu^+}(q)}{\alpha_\theta(q+c+\alpha_\theta)(q-\mu)}
+\frac{\widehat{\nu^+}(q)}{\alpha_\theta(q+c-\alpha_\theta)(q-\mu)}\nonumber\\[1ex]
&\hphantom{\widehat{F_1}(\theta,\mu,q)=}+
\frac{2 \;\widehat{\nu^+}(\alpha_\theta-c)}{(\mu +c-\alpha_\theta)(q+c-\alpha_\theta)(q+c+\alpha_\theta)}
\nonumber\\[1ex]
&\hphantom{\widehat{F_1}(\theta,\mu,q)=}-
\frac{ 2\;\widehat{\nu^+}(\mu)}{(\mu+c-\alpha_\theta)(q-\mu)(q+c+\alpha_\theta)}
\eqp{.}
\end{align}
Let us introduce~:
\begin{equation}
\label{Ctheta}
C_\theta(q):=(q+c+\alpha_\theta)(q+c-\alpha_\theta)=q^2+2cq-2(\lambda+\theta)\eqp{.}
\end{equation}
By a direct calculation we obtain~:
\begin{equation}
\label{trans-F1}
\widehat{F_1}(\theta,\mu,q)=\frac{2}{C_\theta(q)}
\left(
\frac{\widehat{\nu^+}(q)-\widehat{\nu^+}(\mu)}{q-\mu}-
\frac{\widehat{\nu^+}(\mu)-\widehat{\nu^+}(\alpha_\theta-c)}{\mu+c-\alpha_\theta}\right)\eqp{.}
\end{equation}
\goodbreak

b) Let us now compute \mbox{$\displaystyle \widehat{\Lambda_\theta
F}(\theta,\mu,.)(q)$}.

Setting \mbox{$\displaystyle  b=x-a-y$} in (\ref{def-de-lambda}) yields to~:
\begin{align}
&\Lambda_{\theta} F(\theta,\mu,.)(x)\nonumber\\
&= 
\frac{1}{\alpha_{\theta}}\int_{-\infty}^{+\infty} \nu(dy)
\int_0^{+\infty} e^{-c(x-b-y)}
\left( e^{-\alpha_{\theta}|x-y-b|} - e^{-(x+y+b)\alpha_\theta} 
\right) F(\theta,\mu,b) db \nonumber\\ 
&\phantom{==} 
-\frac{1}{\alpha_{\theta}}
\int_{-\infty}^0 \nu(dy)
\int_0^{-y} e^{-c(x-b-y)}
\left( e^{-\alpha_{\theta}(x-y-b)} - e^{-(x+y+b)\alpha_\theta} 
\right) F(\theta,\mu,b) db\nonumber\\
&=H_1F(\theta,\mu,.)(x)+I(\theta,\mu,x)\nonumber\\
&\phantom{==} +\frac{e^{-(c+\alpha_\theta)x}}{\alpha_\theta}
\left(RF(\theta,\mu,.)(\alpha_\theta-c)-RF(\theta,\mu,.)(-\alpha_\theta-c)\right) \eqp{.}   
\end{align}
with $R$ operator defined by (\ref{def-de-R}) and 
\begin{equation}
H_1F(\theta,\mu,.)(x)=\frac{1}{\alpha_{\theta}}\int_{-\infty}^{+\infty}\!\!\!\!\!\! \nu(dy)
\int_{0}^{+\infty}\!\!\!\!\!\! e^{-c(x-b-y)}
e^{-\alpha_{\theta}|x-y-b|}F(\theta,\mu,b) db\eqp{,}       
\end{equation}
\begin{align}
I(\theta,\mu,x)& = 
-\frac{1}{\alpha_{\theta}} 
\int_{-\infty}^{+\infty} \nu(dy)
\int_0^{+\infty}
e^{-c(x-b-y)} e^{-(x+y+b)\alpha_{\theta}} F(\theta,\mu,b)
db\nonumber\\[1ex]
 &=
\frac{1}{\alpha_{\theta}}\;e^{-(c+\alpha_\theta)x}\;\widehat{\nu}(\alpha_\theta-c)\;
\widehat F(\theta,\mu,\alpha_\theta-c)\eqp{.}\hspace{35mm}
\end{align}
Hence
\begin{align}
\label{trans-laplace-lambda}
 \widehat{\Lambda_\theta F}(\theta,\mu,.)(q)&=
\widehat{H_1 F}(\theta,\mu,.)(q)-
\frac{\widehat{\nu}(\alpha_\theta-c)
\widehat F(\theta,\mu,\alpha_\theta-c)}{\alpha_\theta(q+c+\alpha_\theta)}
\nonumber\\[1ex]
&\hphantom{=}+
\frac{RF(\theta,\mu,.)(\alpha_\theta-c)-RF(\theta,\mu,.)(-\alpha_\theta-c)}
{\alpha_\theta(q+c+\alpha_\theta)}
 \eqp{.}        
\end{align}
By definition of $H_1F$, we have~:
\begin{align*}
 \widehat{H_1 F}(\theta,\mu,.)(q)&=
\frac{1}{\alpha_\theta}
\int_{-\infty}^{+\infty}\!\!\!\!e^{cy}\nu(dy)
\int_0^{+\infty}\!\!\!\!e^{cb}
F(\theta,\mu,b)db 
\int_0^{+\infty}\!\!\!\!e^{-(q+c)x}e^{-\alpha_\theta|x-y-b|}dx
\nonumber\\[1ex]
&=
\frac{1}{\alpha_\theta}
\int_{-\infty}^0\!\!\!\! e^{cy}\nu(dy)
\int_0^{-y}\!\!\!\!e^{cb}F(\theta,\mu,b)db 
\int_0^{+\infty}\!\!\!\!e^{-(q+c)x}e^{-\alpha_\theta(x-y-b)}dx
\nonumber\\[1ex]
&\hphantom{=}
+
\frac{1}{\alpha_\theta}
\int_{-\infty}^{+\infty}\!\!\!\!\!\! e^{cy}\nu(dy)
\int_{0\vee (-y)}^{+\infty}\!\!\!\!\!\!e^{cb}F(\theta,\mu,b)db 
\int_0^{b+y}\!\!\!\!e^{-(q+c)x}e^{\alpha_\theta(x-y-b)}dx\nonumber\\
&\hphantom{=}+
\frac{1}{\alpha_\theta}
\int_{-\infty}^{+\infty}\!\!\!\!\!\! e^{cy}\nu(dy)
\int_{0 \vee (-y)}^{+\infty}\!\!\!\!\!\!e^{cb}F(\theta,\mu,b)db
\int_{b+y}^{+\infty}\!\!\!\!e^{-(q+c)x} e^{-\alpha_\theta(x-y-b)}dx
\eqp{.} 
\end{align*}
The $x$-integrals can be computed~:
\begin{align*}
 \widehat{H_1 F}(\theta,\mu,.)(q)&=
\frac{1}{\alpha_\theta}
\int_{-\infty}^0\!\!\!\!e^{(c+\alpha_\theta)y}\nu(dy)
\int_0^{-y}\!\!\!\!e^{(c+\alpha_\theta)b}
F(\theta,\mu,b)db\;\frac{1}{q+c+\alpha_\theta} \nonumber\\[1ex]
&\hphantom{=}
+
\frac{1}{\alpha_\theta}\!
\int_{-\infty}^{+\infty}\!\!\!\!\!\!\!\! e^{(c-\alpha_\theta)y}\nu(dy)
\int_{0\vee (-y)}^{+\infty}\!\!\!\!\!\!\!\!\!\!\!\!e^{(c-\alpha_\theta)b}F(\theta,\mu,b) 
\;\; \frac{1-e^{-(q+c-\alpha_\theta)(b+y)}}{q+c-\alpha_\theta}db
\nonumber\\[1ex]
&\hphantom{=}
+
\frac{1}{\alpha_\theta}
\int_{-\infty}^{+\infty}\!\!\!\!\!\! e^{(c+\alpha_\theta)y}\nu(dy)
\int_{0 \vee (-y)}^{+\infty}\!\!\!\!\!\!\!e^{(c+\alpha_\theta)b}F(\theta,\mu,b)\;\; \frac{e^{-(q+c+\alpha_\theta)(b+y)}}{q+c+\alpha_\theta} db
\eqp{.} 
\end{align*}
By (\ref{def-de-R}) we obtain~:
\begin{align*}
 \widehat{H_1 F}(\theta,\mu,.)(q)&=
 \frac{1}{\alpha_\theta(q+c+\alpha_\theta)}
\left[RF\!(\theta,\mu,.)(\!-\alpha_\theta\!-\!c)+\int_{-\infty}^0\!\!\!\!\!\!\nu(dy)\int_0^{-y}\!\!\!\!\!\!\!\!F(\theta,\mu,b)db\right]
 \nonumber\\[1ex]
&\hphantom{=}
+
\frac{1}{\alpha_\theta(q+c-\alpha_\theta)}\!
\int_{-\infty}^{+\infty}\!\!\!\!\! e^{(c-\alpha_\theta)y}\nu(dy)
\int_{0\vee (-y)}^{+\infty}\!\!\!\!\!e^{(c-\alpha_\theta)b}F(\theta,\mu,b)db \nonumber\\[1ex]
&\hphantom{=}
-
\frac{1}{\alpha_\theta(q+c-\alpha_\theta)}
\int_{-\infty}^{+\infty} e^{-qy}\nu(dy)
\int_{0 \vee (-y)}^{+\infty}e^{-qb}F(\theta,\mu,b)db\nonumber\\
&\hphantom{=}
+
\frac{1}{\alpha_\theta(q+c+\alpha_\theta)}
\int_{-\infty}^{+\infty}e^{-qy} \nu(dy)
\int_{0 \vee (-y)}^{+\infty}e^{-qb}F(\theta,\mu,b)db
\eqp{.} 
\end{align*}
Since
\begin{equation}
\frac{1}{\alpha_\theta(q+c+\alpha_\theta)}-\frac{1}{\alpha_\theta(q+c-\alpha_\theta)}=-\frac{2}{C_\theta(q)}\eqp {,}
\end{equation}
we get~:
\begin{align}
\label{trans-laplace-H1}
 \widehat{H_1 F}(\theta,\mu,.)(q)&=
 \frac{RF(\theta,\mu,.)(-\alpha_\theta-c)}{\alpha_\theta(q+c+\alpha_\theta)} \nonumber\\[1ex]
&\hphantom{=}
+
\frac{\widehat \nu(\alpha_\theta-c)\widehat F(\theta,\mu,\alpha_\theta-c)
-RF(\theta,\mu,.)(\alpha_\theta-c)}{\alpha_\theta(q+c-\alpha_\theta)}
\nonumber\\[1ex]
&\hphantom{=}
-
\frac{2}{C_\theta(q)}
\left(\widehat \nu(q) \widehat F(\theta,\mu,q)-RF(\theta,\mu,.)(q)\right)
\eqp{.} 
\end{align}
Using (\ref{trans-laplace-F}), (\ref{trans-laplace-lambda}) gives~:
\begin{align}
&\widehat F(\theta,\mu,q)
\left(1+\frac{2\widehat \nu (q)}{C_\theta(q)}\right)
= \widehat F_0(\theta,\mu,q)+\widehat F_1(\theta,\mu,q)
+
\frac{2}{C_\theta(q)} RF(\theta,\mu,.)(q)\nonumber\\[1ex]
& \hphantom{=====}+
\frac{2}{C_\theta(q)}
\left(
\widehat \nu(\alpha_\theta-c)\widehat F(\theta,\mu,\alpha_\theta-c)
-RF(\theta,\mu,.)(\alpha_\theta-c) \right)
\eqp{.} 
\end{align}

As  \mbox{$\displaystyle\; C_\theta(q)+2 \widehat \nu(q) = 2 (\varphi(q)-\theta)$}, it is easy to check~:
\begin{align}
(\varphi(q)-\theta) \widehat F(\theta,\mu,q)=&
\frac{1}{2}C_\theta(q)
\left[
\widehat F_0(\theta,\mu,q)+\widehat F_1(\theta,\mu,q)
\right]+
RF(\theta,\mu,.)(q)
\nonumber\\[1ex]
&
+
\widehat \nu(\alpha_\theta-c)\widehat F(\theta,\mu,\alpha_\theta-c)
-RF(\theta,\mu,.)(\alpha_\theta-c)\vphantom{\frac{C_\theta(q)}{2}}
\label{equa-version-calcul} 
\eqp{.} 
\end{align}
Assumptions (\ref{int-sur-R-}) and (\ref{renouv-existence-de-kappa}) imply  the existence of $\gamma^*_0(\theta)$ in $]0,r_\nu[$.
Therefore taking  $q=\gamma^*_0(\theta)$ in (\ref{equa-version-calcul})  brings to~:
\begin{align}
\label{rel-entre-f-chap-et-cteta}
&\widehat \nu(\alpha_\theta-c)\widehat F(\theta,\mu,\alpha_\theta-c)
-RF(\theta,\mu,.)(\alpha_\theta-c)\nonumber\\[1ex]
&=
-\frac{1}{2}C_\theta(\gamma^*_0(\theta))
\left(
\widehat F_0(\theta,\mu,\gamma^*_0(\theta))+\widehat F_1(\theta,\mu,\gamma^*_0(\theta))
\right)- RF(\theta,\mu,.)(\gamma^*_0(\theta))
\eqp{.} 
\end{align}
Determining $\widehat F_0(\theta,\mu,\gamma^*_0(\theta))$ and $\widehat F_1(\theta,\mu,\gamma^*_0(\theta))$ by (\ref{trans-laplace-F0}) et (\ref{trans-F1}), relation (\ref{rel-entre-f-chap-et-cteta}) and (\ref{equa-version-calcul}) imply directly (\ref{equa-Laplace}).
\paragraph{Step 2}\mbox{}
Let $\nu_n$ be the finite measure on $\mathbb R$~:
\begin{equation}
\nu_n(dy):=\nu_{\arrowvert\, {\scriptscriptstyle \mathbb R 
\diagdown
\,  ]-
\frac{1}{n}, \frac{1}{n}[}} (dy) \qquad \forall n \geq 1
\eqp{.}
\end{equation}
We set \mbox{$\displaystyle \lambda_n:=\nu_n(\mathbb R)$}.
We consider $(J_t^n\,,\;t\geq0)$ a compound Poisson process with Lévy measure $\nu_n$, and for any  $n\geq 1$, $x\geq0$ and  $t\geq 0$~:
\begin{align}
X_t^n &=B_t-c_0t+J_t^n  \\[1ex]
T_x^n &=\inf {\{t\geq0\,,\; X_t^n>x\}} \\[1ex]
K_x^n &=X_{T_x^n}^n-x 
\eqp{.}
\end{align}
Let $F_n$ be the Laplace transform of $(T_x^n,K_x^n)$~:
\begin{equation}
F_n(\theta,\mu,x) =
\mathbb E \left(e^{-\theta T_x^n-\mu K_x^n} \un _{\{T_x^n<+\infty\}}\right)\eqp{.}
\end{equation}
By (\ref{equa-Laplace}), the Laplace transform  $\widehat {F_n}(\theta,\mu,.)$ of $ F_n(\theta,\mu,.)$ verifies, for any  $q \in D_0$~:
\begin{align}
\label{transf-Fn}
\widehat{F_n} (\theta,\mu,q)=&
\frac{1}{\varphi_n (q)-\theta}
\left[
\frac{q-\gamma^{*n}_0(\theta)}{2}+
\int_0^{+\infty}\left[\frac{e^{-qy}-e^{-\mu y}}{q-\mu}
\right.\right. \nonumber\\[1ex]
&\left.
-\frac{e^{-\gamma^{*n}_0(\theta)y}-e^{-\mu y}}{\gamma^{*n}_0(\theta)-\mu}\right] \nu_n(dy)\nonumber\\
&\left.
+ R_nF_n(\theta,\mu,.)(q)-R_nF_n(\theta,\mu,.)\!(\gamma^{*n}_0(\theta))\! 
\vphantom{\int}\vphantom{\int} \right] 
\end{align}    
where  $\varphi_n$, $R_n$ and  $\gamma_0^{*n}(\theta)$ are associated with $\nu_n$.
It is well known~:
\begin{alignat*}{2}
&\lim_{n \rightarrow +\infty} T_x^n=T_x \qquad \rm{p.s.}\;\eqp{;}
 && \qquad
\lim_{n \rightarrow +\infty} K_x^n=K_x \qquad \rm{p.s.}\eqp{.}
\end{alignat*}
Consequently~:
\begin{equation*}
\label{limite-Fn}
\lim_{n\rightarrow+\infty} F_n(\theta,\mu,x)= F(\theta,\mu,x) \qquad
\forall \theta, \mu, x \geq 0 \eqp{.}
\end{equation*}
It is easy to check that \mbox{$\displaystyle \lim_{n\rightarrow +\infty} \varphi_n(q)=\varphi(q)$} and \mbox{$\displaystyle \lim_{n\rightarrow +\infty} \gamma_0^{*n}(\theta)=\gamma^*_0(\theta)$}, the proof is left to the reader. Taking the limit $n\rightarrow +\infty$ in (\ref{transf-Fn}), we obtain (\ref{equa-Laplace}).
\finpreuve

\subsection{The particular cases $\nu(]-\infty,0])=0$ and \mbox{$\nu([0,+\infty[)=0$ }}
\label{sec:cas-particuliers}

\noindent {\bf a)}  Let us start with the case $\nu(]-\infty,0])=0$.

\begin{prop}\mbox{} 
\label{prop:q-gamma0}
Assume $\nu_{]-\infty,0[}=0$. Under  (\ref{int-sur-R-}) and (\ref{renouv-existence-de-kappa}), then for any $\mu\geq0$, $\theta \in [0,\, \kappa]$,
\begin{equation}
\label{trans-equi}
\widehat F(\theta,\mu,\rho,q-\gamma_0(\theta))\sim _{q \rightarrow 0}\frac{C_0(\theta,\mu,\rho)}{q}
\end{equation}
where
\begin{enumerate}
\item If $\theta>0$ or if  $\theta=0$ et $ \mathbb E(X_1)\neq 0$,
\begin{align}
\label{def-c-teta-mu}
C_0(\theta,\mu,\rho)=&
\frac{1}{\varphi'(-\gamma_0(\theta))}\!
\left[\frac{-\gamma_0(\theta){-}\gamma^*_0(\theta)}{2}
+\!\!
\int_0^{+\infty}\left[
\frac{e^{(\gamma_0(\theta)-\rho)y}{-}e^{-\mu y}}{-\gamma_0(\theta){+}\rho{-}\mu}\right.\right.\nonumber\\[1ex]
&\hphantom{====}\left.\left.
-\;
\frac{e^{-(\gamma^*_0(\theta)+\rho)y}-e^{-\mu y}}{\gamma^*_0(\theta)+\rho-\mu}\right]\nu(dy)
\right]\eqp{.}
\end{align}
\item If  $\theta=0$ et $\mathbb E(X_1)=0$, 
\begin{equation}
\label{def-C0-teta-mu}
C_0(0,\mu,\rho)=
\frac{1}{\displaystyle \varphi''(0)} \left( 1-\frac{2}{(\rho-\mu)^2}\int_0^{+\infty}(1-e^{(\rho-\mu )y}+(\rho-\mu) y)\nu(dy)\right)\eqp{.}
\end{equation}
\end{enumerate}
\end{prop}

\begin{rem}\mbox{}
\label{rem:q-gamma0}
{\rm
\label{rem:ci}
\begin{enumerate}
\item
In the  companion paper \cite{roynette-vallois-volpi-2003} we prove~:
\begin{equation}
\lim_{x\rightarrow+\infty}e^{\gamma_0(\theta)x}F(\theta,\mu,\rho,x)=
C_0(\theta,\mu,\rho)
\eqp{.}
\end{equation}
\item 
 Since \mbox{$\varphi'(0)=0$} if $\theta=0$ and $\mathbb E(X_1)=0$,   and $\varphi'(-\gamma_0(\theta))<0$ otherwise, this explains why  $C_0(\theta,\mu,\rho)$ is given by two different expressions, (\ref{def-c-teta-mu}) and (\ref{def-C0-teta-mu}).
\item 
The constant $C_0(\theta,\mu,\rho)$ is positive because  \mbox{$\displaystyle \mu \rightarrow \widehat F(\theta,\mu,q-\!\gamma_0(\theta))$} is decreasing and \mbox{$\displaystyle \lim_{\mu\rightarrow+\infty}C_0(0,\mu,\rho)=\frac{1}{\varphi''(0)}>0$} \;  if $\theta=0$ and \mbox{$\mathbb E(X_1)=0$}, \goodbreak and  $\displaystyle \lim_{\mu\rightarrow+\infty}C_0(\theta,\mu,\rho)=-\;\frac{\gamma_0(\theta)+\gamma_0^*(\theta)}{2\varphi'(-\gamma_0(\theta))}>0$ otherwise.
\item
The constant $C_0(0,0,0)$ can be computed explicitely~:
\begin{equation}
C_0(0,0,0)=-\;\frac{\varphi'(0)}{\varphi'(-\gamma_0(0))}
\end{equation}
when $\mathbb E(X_1)<0$, and $C_0(0,0,0)=1$ otherwise.
\end{enumerate}
}
\end{rem}

\noindent{\bf Proof of Proposition~\ref{prop:q-gamma0}}\mbox{}

\noindent Once more we only deal with $\rho=0$, and $F(\theta,\mu,.)$ stands for $F(\theta,\mu,\rho,.)$.

\noindent $1)$ We first suppose   $\theta> 0$ or $\theta=0$ and $ \mathbb E(X_1)\neq 0$.

\noindent  Recall (cf. Remark~\ref{rem:ci}, point $2$)   
\mbox{$\varphi'(-\gamma_0(\theta))\neq 0$}, then~:
\begin{equation}
\label{equi-phi-prim}
\varphi(q-\gamma_0(\theta))-\theta=\varphi(q-\gamma_0(\theta))-\varphi(-\gamma_0(\theta))\sim_{q \rightarrow0}q\varphi'(-\gamma_0(\theta))\eqp{.}
\end{equation}
Replacing in  (\ref{trans-nu+}) $q$ by $q-\gamma_0(\theta)$ and taking the limit as  $q\rightarrow 0$, we conclude immediately that (\ref{trans-equi}) holds.

\noindent $2)$ If  $\theta=0$ and  $\mathbb E(X_1)=0$, then $\gamma_0(0)=0$, $\varphi'(0)=0$ and 
\begin{equation}
\label{equi-q} 
\varphi(q)\sim_{q \rightarrow 0} \frac{q^2}{2} \varphi''(0)\eqp{.}
\end{equation}
(\ref{trans-equi}) follows easily.\finpreuve

\mbox{}\\
\noindent {\bf b) }  We now briefly investigate the case $\nu(]0,+\infty[)=0$.\\
We observe that $K_x=L_x=0$, then \mbox{$\displaystyle F(\theta,\mu,\rho,.)=F(\theta,0,0,.)$}. We can check that the function 
 \mbox{$\displaystyle G_\theta : x  \rightarrow e^{-\gamma_0(\theta)x}$}  verifies (\ref{equa-Laplace}) and so is the unique solution of the  functional equation (\ref{equa-fonctionnelle}) (For a  proof, see \cite{volpi03}).


\sectio{A new functional equation verified by $F$}
\label{sec:new-funct-equa}


Suppose that $\nu$ satisfied the assumption given in Theorem~\ref{thm:equa-Laplace}. If  $\nu(\mathbb R)<+\infty$, since the   Laplace transformation  is  one-to-one,  (\ref{equa-fonctionnelle}) is equivalent to (\ref{equa-Laplace}). But relation (\ref{equa-Laplace}) remains valid when $\nu(\mathbb R)=+\infty$. This brings us to ask what is the relation involving $F$ induced by (\ref{equa-Laplace}). In other words is it possible to inverse (\ref{equa-Laplace})? That strengthen the role of equation (\ref{equa-Laplace}) and also the approach we have  developed previously via the Laplace transform of $F(\theta,\mu,\rho,.)$.  Let $\widetilde L$ be the operator~:
\begin{equation}
\label{jolie-nouvelle-equa}
\widetilde Lf(x) = \frac{1}{2} f''(x)+cf'(x)+\int_{-\infty}^{+\infty} (f(x-y)-f(x)) \nu(dy)\eqp{.}
\end{equation}
We notice that  $\widetilde L$ is the formal adjoint of the infinitesimal generator $L$ of ($\X$).

\begin{thm}\mbox{}
\label{thm:seconde-equa-fonc}
Suppose that $\nu$ satisfies the hypotheses given in Theorem~\ref{thm:equa-Laplace} and morever $\displaystyle \int_{-1}^1|y|\nu(dy)<\infty$ . Then 
\begin{equation}
\label{seconde-equa-fonc}
\widetilde LF(\theta,\mu,\rho,x)-\theta F(\theta,\mu,\rho,x)=g(\mu,\rho,x)\; ;  \quad x>0 \eqp{.}
\end{equation}
where
\begin{equation}
g(\mu,\rho,x):= -\left(\int_{[x,+\infty[}e^{-\mu(y-x)-\rho x}\nu(dy)\right)\un_{\{x>0\}}\eqp{,}
\end{equation}
 with the boundary conditions~:
\begin{align}
\label{condition-sur-F}
&F(\theta,\mu,\rho,0)=1\eqp{,}\\
&F'(\theta,\mu,\rho,0_+)=-2 \left[c+\frac{\gamma_0^*(\theta)}{2}+\frac{1}{\gamma_0^*(\theta)+\rho-\mu}\int_0^{+\infty}\!\!\!\!(e^{-(\gamma_0^*(\theta)+\rho)y}-e^{-\mu y}) \nu(dy)\right.\nonumber\\
&\hphantom{=========}\left. +\int_{-\infty}^0 \nu(dy) \int_0^{-y} e^{-\gamma_0^*(\theta)(y+b)} F(\theta,\mu,\rho,b)db\right]
\label{derivee-F-0+}
\eqp{.}
\end{align}
\end{thm}

\begin{rem}\mbox{}\\
{\rm
\noindent $1.$\;  In (\ref{seconde-equa-fonc}) and (\ref{derivee-F-0+})  the derivatives are $x$-derivatives.\\
\noindent $2.$\;  Suppose $\nu(]-\infty,0])=0$. Then the last term in the right hand-side of (\ref{derivee-F-0+}) cancels and $F'(\theta,\mu,\rho,0_+)$ only depends on $\nu$. Consequently  (\ref{seconde-equa-fonc}) is a classical integro-differential linear equation. If morever $\mu=\rho=0$, then (\ref{derivee-F-0+}) reduces to $\displaystyle F'(\theta,0,0,0_+)=-\frac{2\theta}{\gamma_0^*(\theta)}$. If additionally $\theta=\mu=\rho=0$, the ruin probability $F(0,0,0,.)$ solves~:
\begin{equation}
\label{sec-equa-nu+}
\widetilde L F(0,0,0,x)=-\nu([x,+\infty[)\;, \qquad x>0 \eqp{.}
\end{equation}
with $F(0,0,0,0)=1$ and $F'(0,0,0,0_+)=0$ if $\mathbb E(X_1)\geq 0$ and 

\noindent \mbox{$F'(0,0,0,0_+)=-2 \varphi'(0)=2\mathbb E(X_1)$} if $\mathbb E(X_1)< 0$.
It is easy to check that $x\rightarrow 1$ is the unique  solution of (\ref{sec-equa-nu+}) with the  boundary conditions (\ref{condition-sur-F}) and  (\ref{derivee-F-0+}) when $\mathbb E(X_1)\geq 0$.\\
\noindent $3.$ \;
Suppose $\nu([0,+\infty[)=0$. Then $g(\mu,\rho,.)$ cancels and (\ref{seconde-equa-fonc}) reduces to~:
\begin{equation}
\label{sec-equa-nu-}
\widetilde LF(\theta,\mu,\rho,x)-\theta F(\theta,\mu,\rho,x)=0\;;\quad x>0\eqp{.}
\end{equation}
It is easy to verify that $x\rightarrow e^{-\gamma_0(\theta)x}$ is the unique solution of (\ref{sec-equa-nu-}), (\ref{condition-sur-F}) and  (\ref{derivee-F-0+}).\\
\noindent $4.$\; 
Obviously (\ref{seconde-equa-fonc}) may be written as~:
\begin{equation}
\label{nouvelle-equa}
\frac{1}{2}F''(\theta,\mu,\rho,x)+c F'(\theta,\mu,\rho,x)-\theta F(\theta,\mu,\rho,x)=h(x)\eqp{,}
\end{equation}
where
\begin{equation}
\label{def-h}
h(x)=g(\mu,\rho,x)-\int_{-\infty}^{+\infty} (F(\theta,\mu ,\rho,x-y)-F(\theta,\mu,\rho,x))\nu(dy)
\end{equation}
Considering (\ref{nouvelle-equa}) as a linear differential equation with given data $h$, and integrating with the method of variation of  parameter, we obtain~:
\begin{align}
\label{nvelle-equa-finale}
F(\theta,\mu,\rho,x)=&\frac{e^{\alpha_1 x}}{\sqrt{c^2+2\theta}}\left[(\alpha_2-F'(\theta,\mu,\rho,0_+))-\int_0^x e^{-\alpha_1 y}h(y) dy\right]\nonumber\\
&
+\frac{e^{\alpha_2 x}}{\sqrt{c^2+2\theta}}\left[(F'(\theta,\mu,\rho,0_+)-\alpha_1)-\int_0^x e^{-\alpha_2 y}h(y) dy\right]
\end{align}
where $\displaystyle \alpha_1=-c+\sqrt{c^2+2\theta}$ and 
$\displaystyle \alpha_2=-c-\sqrt{c^2+2\theta}$.\\
\noindent $5.$\; 
We point out that (\ref{nvelle-equa-finale}) can be written as~:
\begin{equation}
F(\theta,\mu,\rho,x)= F_2(\theta,\mu,\rho,x)+\overline{\Lambda}_\theta F(\theta,\mu,\rho,.)(x)\eqp{,}
\end{equation}
where
\begin{align}
& F_2(\theta,\mu,\rho,x):=\frac{e^{\alpha_1 x}}{\sqrt{c^2+2\theta}}\left[\alpha_2-\int_0^x e^{-\alpha_1 y}g(\mu,\rho,y) dy\right]\nonumber\\[1ex]
&+
\frac{e^{\alpha_2 x}}{\sqrt{c^2+2\theta}}\left[
-\alpha_1-\int_0^x e^{-\alpha_2 y}g(\mu,\rho,y) dy\right]
\frac{\left(2c+\gamma_0^*(\theta)\right)
\left(e^{\alpha_1 x}-e^{\alpha_2 x}\right)}{\sqrt{c^2+2\theta}}\nonumber\\[1ex]
&+
\frac{2\left(e^{\alpha_1 x}-e^{\alpha_2 x}\right)}{\sqrt{c^2+2\theta}
\left(\gamma_0^*(\theta)+\rho-\mu\right)}\int_0^{+\infty}\!\!\!\!(e^{-(\gamma_0^*(\theta)+\rho)y}-e^{-\mu y}) \nu(dy)
\end{align}
and $\overline{\Lambda}_\theta$ is the linear operator~:
\begin{align}
&\overline{\Lambda}_\theta f:= 
\frac{2 \left(e^{\alpha_1 x}-e^{\alpha_2 x}\right)}{\sqrt{c^2+2\theta}}
\int_{-\infty}^0 \nu(dy) \int_0^{-y} e^{-\gamma_0^*(\theta)(y+b)} F(\theta,\mu,\rho,b)db\nonumber\\[1ex]
&+\frac{e^{\alpha_1 x}}{\sqrt{c^2+2\theta}}\int_0^x e^{-\alpha_1 y}
\left(\int_{-\infty}^{+\infty} (F(\theta,\mu ,\rho,y-z)-F(\theta,\mu,\rho,y))\nu(dz)\right)
dy\nonumber\\[1ex]
&+
\frac{e^{\alpha_2 x}}{\sqrt{c^2+2\theta}}
\int_0^x e^{-\alpha_2 y}
\left(\int_{-\infty}^{+\infty} (F(\theta,\mu ,\rho,y-z)-F(\theta,\mu,\rho,y))\nu(dz)\right)
 dy
\end{align}
}
\end{rem}

\noindent{\bf Proof of Theorem~\ref{thm:seconde-equa-fonc}}\mbox{}

Multiplying both sides of (\ref{equa-Laplace}) by $\varphi(q)-\theta$, we obtain~:
\begin{align}
\label{multiplier-par-phi}
(\varphi(q)-\theta)&\widehat F(\theta,\mu,\rho,q)-\frac{q}{2}+ \frac{\gamma_0^*(\theta)}{2} =\nonumber\\
&\int_0^{+\infty}\left[\frac{(e^{-(q+\rho)y}-e^{-\mu y})}{q+\rho-\mu}
-\frac{(e^{-(\gamma_0^*(\theta)+\rho)y}-e^{-\mu y})}{\gamma_0^*(\theta)+\rho-\mu}\right]\nu(dy)\nonumber\\
&+RF(\theta,\mu,\rho,.)(q)-RF(\theta,\mu,\rho,.)(\gamma_0^*(\theta))\eqp{.}
\end{align}
We observe that the left hand-side and the right hand-side of 
(\ref{multiplier-par-phi}) are Laplace transforms. This leads to (\ref{seconde-equa-fonc}), the details are left to the reader.\finpreuve

\appendix


\sectio{Annex~: Properties of  $\varphi$ and  $\varphi_\theta$}

\label{A1-fonction-phi}



\subsection{Properties of  $\varphi$ }

\label{phi-etude}


Recall that~:

\begin{equation}
\label{A-def-de-phi}
\varphi(q)=\frac{q^2}{2}+cq+\int_{-\infty}^{+\infty}\left(e^{-qy}-1+qy\un_{\{|y|\leq 1\}}\right) \nu(dy)\eqp{.}
\end{equation}

 (\ref{A-def-de-phi}) implies that  $\varphi(q)$  exists if~:

\begin{equation}
\int_\mathbb R \un_{\{|y|\geq1\}}|e^{-qy}|\nu(dy)<+\infty,\quad q\in\mathbb C \eqp{.}
\end{equation}

Recall that $r_\nu$ is defined by (\ref{def-de-r}). Let

\begin{equation}
\label{A-def-de-r*}
r^*_\nu := \sup \left\{s \geq 0 ;\; \int_{-\infty}^{-1} e^{-sy}
\nu (dy )\;<\;+\infty \right\}\eqp{,}
\end{equation}

From now on,  we suppose~:

\begin{equation}
\label{A-hyp-sur-r-et-r*}
r_\nu > 0 \quad {\rm et}\quad r_\nu^*>0\quad (r_\nu\;{\rm or }\;r_\nu^* \;{\rm may\;  be \;infinite}) \eqp{.}
\end{equation}

By definition  $\varphi$ is a convex function defined on $]-r_\nu,r^*_\nu[$.

\noindent Morever

\begin{equation}
\varphi'(0)=\mathbb E(X_1)\eqp{.}
\end{equation}

If  $r_\nu<+\infty$ or  $r_\nu^*<+\infty$,  we extend  $\varphi$ as follows~:

\begin{equation}
 \varphi(-r_\nu)=\lim_{q\rightarrow-r_\nu}\varphi(q)\quad{\rm et}\quad 
\varphi(r_\nu^*)=\lim_{q\rightarrow r_\nu^*}\varphi(q)\eqp{.}
\end{equation}

We plot below (see figures~\ref{fig:profit-net}, \ref{fig:non-profit-net} and  ~\ref{fig:egalite}) the graph of $\varphi$, distinguishing three cases~: 

\mbox{$\displaystyle \mathbb E(X_1)<0$}, \mbox{$\displaystyle \mathbb E(X_1)>0$} and \mbox{$\displaystyle \mathbb E(X_1)=0$}.

\begin{figure}[h]

\begin{minipage}{0.3\linewidth}
\centering
\scalebox{0.7}{\input{profit-net.pstex_t}}
\subfigure{a) $-r_\nu<-\gamma_0(0)$}
\end{minipage}\goodgap
\begin{minipage}{0.3\linewidth}
\centering
\scalebox{0.7}{\input{pn-rnu-egal.pstex_t}}
\subfigure{b) $-\gamma_0(0)=-r_\nu$}
\end{minipage}
\begin{minipage}{0.3\linewidth}
\centering
\scalebox{0.7}{\input{non-gamma0-pn.pstex_t}}
\subfigure{c) one zero only}
\end{minipage}
\caption{Graph of $\varphi$, $\mathbb E(X_1)<0$}
\label{fig:profit-net}
\end{figure}

\begin{figure}[h]
\begin{minipage}{0.3\linewidth}
\centering
\scalebox{0.7}{\input{non-pn.pstex_t}}
\subfigure{a) $\gamma_0^*(0)<r_\nu^*$}
\end{minipage}\goodgap
\begin{minipage}{0.3\linewidth}
\centering
\scalebox{0.7}
{\input{non-pn-rnux-egal.pstex_t}}
\subfigure{b) $\gamma_0^*(0)=r_\nu^*$}
\end{minipage}
\begin{minipage}{0.3\linewidth}
\centering
\scalebox{0.7}
{\input{non-gamma-non-pn.pstex_t}}
\subfigure{c) one zero only}
\end{minipage}
\caption{Graph of $\varphi$, $\mathbb E(X_1)>0$}
\label{fig:non-profit-net}
\end{figure} 
\begin{figure}[h]
\centering
\scalebox{0.7}
{\input{egalite.pstex_t}}
\caption{Graph of $\varphi$, $\mathbb E(X_1)=0$}
\label{fig:egalite}
\end{figure}

\noindent When $\varphi$ has two zeros (may be a double zero) in $[-r_\nu, r_\nu^*]$, $-\gamma_0(0)$ (resp. $\gamma_0^*(0)$) will denote the smallest (resp. biggest) one.

\begin{prop}\mbox{}
\label{prop:etude-de-phi}
\begin{enumerate}
\item
In particular in cases Fig~\ref{fig:profit-net} a, Fig~\ref{fig:non-profit-net} a  and Fig~\ref{fig:egalite}, we have~:
\begin{equation}
-r_\nu<-\gamma_0(0)\leq 0\leq \gamma_0^*(0)<r_\nu^*\eqp{.}
\end{equation}
\item The set $\displaystyle \{s\in [-r_\nu,0]\;/\; \varphi(s)<0\}$ is an  interval,  being non empty as soon as   $\mathbb E(X_1)<0$. 
\end{enumerate}
\end{prop}
%


\subsection{Properties of  $\varphi_\theta:q\rightarrow\varphi(q)-\theta$ }
\label{phi-theta-etude}

\begin{figure}[h]
\centering
\scalebox{0.6}
{\input{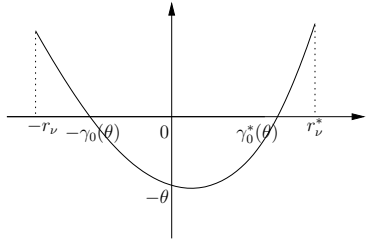}}
\caption{Graph of $\varphi_\theta=\varphi-\theta$}
\label{fig:theta}
\end{figure}

Let us briefly mention properties of $\varphi_\theta=\varphi-\theta$, $\theta>0$. For simplicity we restrict ourselves  to  Fig~\ref{fig:profit-net} a, Fig~\ref{fig:non-profit-net} a  and Fig~\ref{fig:egalite}. The graph of $\varphi_\theta$ is given by Fig~\ref{fig:theta}.

Then there exists $\kappa>0$, such that for any $\theta \in\; ]0,\kappa]$, $\varphi_\theta$ has an unique positive ( resp. negative) zero denoted $\gamma_0^*(\theta)$ ( resp.  $-\gamma_0(\theta)$) and~: 
\begin{equation}
-r_\nu<-\gamma_0(\theta)<-\gamma_0(0)\leq 0 \leq \gamma_0(0)<\gamma_0^*(\theta)<r_\nu^*\eqp{.}
\end{equation}


\subsection{The zeros of  $C_\theta (q)$}
\label{polynome-C-theta}


\noindent Assume that $\lambda=\displaystyle \nu(\mathbb R)<+\infty $.\\

We notice that for any $q\in]-r_\nu,r_\nu^*[$, we have~:
\begin{equation}
\varphi_\theta(q)=\frac{1}{2}C_\theta(q)+\widehat \nu(q)\eqp{,}
\end{equation}
where $\widehat \nu$ is the Laplace transform of $\nu$, i.e.~:
\begin{equation}
\widehat \nu(q):=\int_{-\infty}^{+\infty}e^{-qy}\nu(dy)\eqp{,}
\end{equation}
and   $C_\theta(q)$ is the  polynomial function~:
\begin{equation}
C_\theta (q)=q^2+2c_0q-2(\lambda+\theta)\eqp{,}
\end{equation}
with $\displaystyle c_0=c+\int_\mathbb R y \un_{\{|y|\leq 1\}}\nu(dy)$.
Then
\begin{equation}
\label{phi-theta-et-C-theta}
\varphi_\theta(q)>\frac{1}{2}C_\theta(q),\quad q\in]-r_\nu,r_\nu^*[\eqp{.}
\end{equation}
Recall
\begin{equation}
\label{factorisation-de-C-theta}
C_\theta(q)=(q+c_0+\alpha_\theta)(q+c_0-\alpha_\theta)\eqp{,}
\end{equation}
where \mbox{$\displaystyle \alpha_\theta=\sqrt{c_0^2+2(\lambda+\theta)}$}.\\
Let us  summarize the results in the following.

\begin{prop}
\label{prop:zeros-de-C-theta}
\quad 
Suppose  $\theta \geq 0$. Then~:
\begin{itemize}
 \item[$(i)$]   the two real  zeros of $C_\theta(q)$  are  \mbox{$-\alpha_\theta-c_0<0$} and  \mbox{$\alpha_\theta-c_0>0$},
\item[$(ii)$]
the set $\displaystyle \{s\in\;]-r_\nu,r_\nu^*[\;/\;\varphi_\theta(s)<0\}$ is an  interval included in 

\noindent \mbox{$]-c_0-\alpha_\theta,\alpha_\theta-c_0[$},
\item[$(iii)$]

if $\varphi_\theta$ has two zeros $-\gamma_0(\theta)$ et $\gamma_0^*(\theta)$ in \mbox{$[r_\nu,r_\nu^*]$} then~:
\begin{equation}
\label{comparaison-de-racines}
-c_0-\alpha_\theta<-\gamma_0(\theta)\leq 0\leq \gamma_0^*(\theta)<\alpha_\theta-c_0\eqp{.}
\end{equation}
\end{itemize}
\end{prop}

\begin{figure}[t]
\centering
\scalebox{0.9} 
{\input{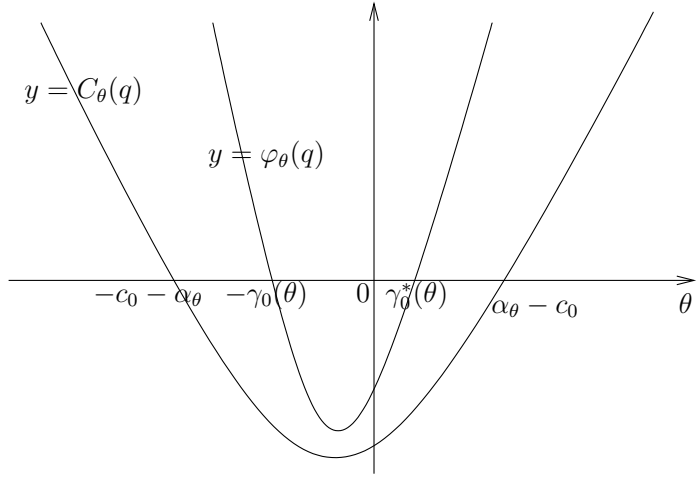}}
\caption{Comparison of the zeros of   $C_\theta(q)$  and those of   $\varphi_\theta$}
\label{fig:zeros-de-C(q)}
\end{figure}

\clearpage
\bibliographystyle{plain}
\bibliography{Biblio}

\begin{thebibliography}{10}

\bibitem{Bertoin96}
J.~Bertoin.
\newblock {\em L{\'e}vy process}.
\newblock Cambridge Univ. Press, 1996.
\newblock Cambridge tracts in mathematics vol. 121.

\bibitem{BertoinDoney94}
J.~Bertoin and R.A. Doney.
\newblock Cram\'er's estimate for {L}évy processes.
\newblock {\em Statistics and Probability Letters}, 21:363--365, 1994.

\bibitem{Cramer30}
H.~Cramér.
\newblock {\em On the mathematical Theory of Risk}.
\newblock Skandia Jubilee Volume, Stockholm, 1930.

\bibitem{Cramer55}
H.~Cramér.
\newblock {\em Collective Risk Theory}.
\newblock Skandia Jubilee Volume, Stockholm, 1955.

\bibitem{DassiosEmbrechts89}
A.~Dassios and P.~Embrechts.
\newblock Martingales in insurance risk.
\newblock {\em Commun. Statistics Stochastics Models}, 5(2):181--217, 1989.

\bibitem{DelbaenHaezendonck85}
F.~Delbaen and J.~Haezendonck.
\newblock Inversed martingales in risk theory.
\newblock {\em Insurance: Mathematics and Economics}, 4:201--206, 1995.

\bibitem{DicksonWaters92}
D.C.M. Dickson and H.R. Waters.
\newblock The probability and severity of ruin in finite and infinite time.
\newblock {\em ASTIN Bulletin}, 22:177--190, 1992.

\bibitem{DozziVallois97}
M.~Dozzi and P.~Vallois.
\newblock Level crossing times for certain processus without positive jumps.
\newblock {\em Bulletin des Sciences Mathématiques}, 121:355--376, 1997.

\bibitem{DufresneGerber90}
F.~Dufresne and H.U. Gerber.
\newblock Risk theory for the compound {P}oisson process that is perturbed by
  diffusion.
\newblock {\em Insurance: Mathematics and Economics}, 10(1991):51--59, 1990.

\bibitem{DufresneGerberShiu91}
F.~Dufresne, H.U. Gerber, and E.~Shiu.
\newblock Risk theory with the gamma process.
\newblock {\em ASTIN Bulletin}, 21:177--192, 1991.

\bibitem{Feller71}
W.~Feller.
\newblock {\em An Introduction to Probability Theory and its Applications},
  volume~II.
\newblock John Wiley and Sons, New York, 2 edition, 1971.

\bibitem{Gerber70}
H.U. Gerber.
\newblock An extension of the renewal equation and its application in the
  collective theory of risk.
\newblock {\em Skandinavisk Aktuarietidskrift}, 53:205--210, 1970.

\bibitem{Gerber73}
H.U. Gerber.
\newblock Martingales in risk theory.
\newblock {\em Mitteilungen der Schweizerischen Vereinigung der
  Versicherungsmathematiker}, pages 205--216, 1973.

\bibitem{GerberShiu97}
H.U. Gerber and E.S.W. Shiu.
\newblock The joint distribution of the time of ruin, the surplus immediately
  before ruin and the deficit at ruin.
\newblock {\em Insurance: mathematics and Economics}, 21:129--137, 1997.

\bibitem{gradshteyn80}
I.S. Gradshteyn and I.M. Ryzhik.
\newblock {\em Table of integrals, series and products}.
\newblock Academic Press, INC, Orlando, corrected and enlarged edition, 1980.

\bibitem{HuzakPerman03}
M.~Huzak, M.~Perman, H.~Sikic, and Z.~Vondracek.
\newblock {R}uin probabilities and decompositions for general perturbed risk
  processes.
\newblock preprint, 2003.

\bibitem{Karatszas91}
I.~Karatszas and S.E. Shreve.
\newblock {\em Brownian motion and Stochastic Calculus}.
\newblock Springer-Verlag, New York, second edition, 1991.
\newblock Graduate Texts in Mathematics, 113.

\bibitem{KouWang01}
S.G. Kou and H.~Wang.
\newblock First passage times of a jump diffusion process.
\newblock {\em Preprint}, 2001.

\bibitem{Lundberg03}
F.~Lundberg.
\newblock {\em I- Approximerad Framställning av Sannolikhetsfunktionen. II-
  Aterförsäkering av Kollectivrisker}.
\newblock Almqvist and Wiksell, Uppsala, 1903.

\bibitem{NguyenYor01}
L.~Nguyen and M.~Yor.
\newblock Wiener {H}opf factorization and the pricing of barrier of and
  lookback options under general {L}évy processes.
\newblock {\em Preprint}, 2001.

\bibitem{Norberg}
R.~Norberg.
\newblock Ruin problems with assets and liabilities of diffusion type.
\newblock preprint.

\bibitem{RolskiSchmidli99}
T.~Rolski, H.~Schmidli, V.~Schmidt, and J.~Teugels.
\newblock {\em Stochastic Processes for Insurance and Finance}.
\newblock Wiley Series in Probability and Statistics. John Wiley and Sons,
  Chichester, New York, Weinheim, Brisbane, Singapore, Toronto, 1999.

\bibitem{roynette-vallois-volpi-2003}
B.~Roynette, P.~Vallois, and A.~Volpi.
\newblock {L}évy processes~: {H}itting time, overshoot and undershoot; {II}-
  {A}symptotic behavior.
\newblock {\em Preprint}, 2003.

\bibitem{Taylor76}
G.C. Taylor.
\newblock Use of differential and integral inequalities to bound ruin and
  queueing probabilities.
\newblock {\em Scandinavian Actuarial Journal}, pages 57--76, 1976.

\bibitem{volpi03}
A.~Volpi.
\newblock {\em Processus associés à l'équation de diffusion rapide; Etude
  asymptotique du temps de ruine et de l'overshoot}.
\newblock Univ. Henri Poincaré, Nancy I, Vandoeuvre les Nancy, 2003.
\newblock Thèse.

\end{thebibliography}

\end{document}